\documentclass[10pt]{article}
\usepackage[colorlinks=true, pdfstartview=FitV, linkcolor=blue, citecolor=blue, urlcolor=blue]{hyperref}
\usepackage{amsmath,amsfonts,latexsym,amssymb}
\usepackage[latin1]{inputenc}
\newtheorem{theorem}{Theorem}[subsection]
\newtheorem{lemma}[theorem]{Lemma}
\newtheorem{proposition}[theorem]{Proposition}
\newtheorem{coro}[theorem]{Corollary}
\newtheorem{conj}[theorem]{Conjecture}
\newcommand{\grf}{\pi_1(S)}
\newcommand{\proof}{{\sc Proof : }}
\newcommand{\qed}{{\sc Q.e.d.}}
\newcommand{\rks}{{\sc Remarks}}
\newcommand{\auteur}{
\vskip 2truecm
\centerline{Univ. Paris-Sud,}
\centerline{Laboratoire de Mathématiques, Orsay F-91405 Cedex }
\centerline{CNRS, Orsay cedex, F-91405}}
\begin{document}
\title{Flat Projective Structures on Surfaces and \\  Cubic Holomorphic Differentials}
\author{François LABOURIE \thanks{Univ. Paris-Sud, Laboratoire de Mathématiques, Orsay F-91405 Cedex; CNRS, Orsay cedex, F-91405}}
\maketitle

\newcommand\inn{_{n\in\mathbb N}}
\newcommand\ci{C^\infty}
\newcommand\g{\gamma}
\newcommand{\nabar}{\overline\nabla}
\newcommand\f{\phi}
\newcommand\la{\lambda}
\newcommand\cG{{\mathcal G}}
\newcommand\D{\Delta}
\newcommand\pr{\partial}
\newcommand\e{\epsilon}
\newcommand\cF{{\mathcal F}}
\renewcommand\a{\alpha}
\newcommand\RR{\mathbb R}
\newcommand\PP{\mathbb P}
\newcommand\cU{\cal U}
\newcommand\cC{\cal C}
\newcommand\tr{{\rm trace}}
\newcommand\CI{C^{\infty}}
\renewcommand\sp{flat projective structure}
\newcommand\rpn{{\mathbb {RP}}^{n}}
\newcommand\rpd{{\mathbb {RP}}^{2}}
\newcommand\wtm{\widetilde M}
\newcommand\wtn{\widetilde \nabla}
\newcommand\lbreak{\medskip }
\newcommand\lb{\medskip }
\newcommand{\vecteur}[1]{
\left( \begin{array}{c} 
#1 \end{array}\right)
}
\newcommand{\mapping}[4]
{
\left\{
\begin{array}{rcl}
#1 &\rightarrow& #2\\
#3 &\mapsto& #4 
\end{array}
\right.
}
\newcommand{\dt}[1]{\frac{d#1}{dt}}

\section{Introduction}
The purpose of this article is to give an interpretation of real projective structures and associated cohomology classes in terms of connections, sections, etc.  satisfying  elliptic partial differential equations in the  spirit of Hodge theory. We shall also give an application of these results as the uniqueness of a minimal surface in a symmetric space.

We recall briefly that a flat real projective structure on a surface $S$ is an atlas with values in $\mathbb{RP}^2$ and  coordinates changes in  $SL(3,\mathbb R)$. Associate to such a structure is a a {\em holonomy representation} of $\grf$ with values in $SL(3,\mathbb R)$, and {\em a developing map}, defined  from $\tilde S$, the universal cover of $S$, with values in $\mathbb{RP}^2$ and equivariant under $\rho$. Finally, the structure is said to be {\em convex} if the image of the developing map is  a convex set. 

Convex projective structures have been extensively studied by Choi Suhyoung and William Goldman in \cite{WGo} and  \cite{Ch-Go} for instance. We summarise some of their major results in the following

\begin{theorem}{\sc[Choi-Goldman]}
Every convex structure on $S$ is determined by its holonomy representation. Moreover, if a representation of $\grf$ in $SL(3,\mathbb R)$ can be deformed into a discrete faithful representation in $SO(2,1)$, then  it is the holonomy of a convex structure on $S$.
\end{theorem}
In the introductory Section \ref{sec:intro}, we give various points of view on projective structures.

Most of the results of this paper can be stated as bijections between moduli spaces and set of solutions of certain equations. Of course the important point, not always clear in the statements, is the construction of the bijection. We now give a sketch of the content of this article. Most of the material of this article is new, although some results were announced a long time ago

\subparagraph{Convex projective structures and cubic holomorphic differentials.}
Concerning convex projective structures, we prove

\begin{theorem}\label{intro:A}
There exists a mapping class group equivariant homeomorphism between the moduli space of convex structures on $S$ and the moduli space of pairs $(J,Q)$ where $J$ is a complex structure on $S$ and $Q$ is a cubic holomorphic differential on $S$ with respect to $J$.
\end{theorem}
This result is a combination of Theorem \ref{rp->cub} and \ref{cub->rp}. This result was announced in \cite{FL2}.  The proof presented in the present paper uses a self contained approach. However, this result, as it is  explained in Section \ref{sec:aff}, can be obtained as a consequence of difficult results in affine differential geometry of Cheng and Yau \cite{CY1} and  \cite{CY2} later completed and clarified by the work of Gigena \cite{Gi}, Sasaki \cite{Sa} and A.M.\ Li \cite{Li}, \cite{Li2}. In \cite{JL}, John Loftin also proves and extends this result using this affine differential geometric interpretation.

\subparagraph{Projective structure and cohomology classes.}

In Section \ref{secOp}, we associate to every flat projective structure of holonomy $\rho$, a non empty cone in $H^1_\rho(\mathbb R^3)$. We also prove the this cone helps to distinguish convex structures from others. Indeed, this cone contains $0$ if and only if the structure is convex.

\subparagraph{Cohomology classes and complex structures.}

We obtain two results which parametrise the moduli space of representations of the surface group in the affine space in dimension 3: Theorem \ref{pcohom} which defines a map from $H^1_\rho(E)$ to Teichmüller space, and Theorem \ref{phodge} which can be though as a generalisation of the Eichler-Shimura isomorphism in this context.

\subparagraph{Dualities and symmetries.}

We show that the above results, more specially Theorem \ref{pcohom} and Theorem \ref{phodge}, give rise to unexpected symmetries of  the moduli space of representations of the surface group in the affine space in dimension 3.

\subparagraph{Higgs bundle interpretation.}
We also interpret Theorem \ref{intro:A} in the context of Higgs bundle theory. For instance, using this result, we obtain that the energy map on Teichmüller space associated to the holonomy representation of a projective structure has a unique critical point which is minimum. We note that this function is proper according to \cite{FLEnergy}. We obtain in particularly Corollary \ref{cor:min} that states the existence and uniqueness of a minimal surface in some symmetric spaces.

\subparagraph{Holomorphic interpretation.}  In Section \ref{sec:holo}, we explain how Theorem \ref{intro:A} can be interpreted as the existence an uniqueness of an equivariant "holomorphic"  curve in $SL(3,\mathbb R)/SL(2,\mathbb R)$.

\subparagraph{Appendices.} Finally, Appendices \ref{lapla} and \ref{ma} contain  compactness results for partial differentials equations appearing in the paper, arising as  consequences of an holomorphic interpretations in the spirit of \cite{maudin}, and \cite{labgafa}.

\tableofcontents
\section{Projective structures}\label{sec:intro}

\subsection{First definitions.}\label{def}

Let $M$ be a $n$-dimensional manifold.

\subsubsection{Projectively equivalent connections}
On a manifold $M$ two connections are said to be {\em projectively equivalent} if they have the same geodesics, up to  parametrisations.  In dimension greater than 2, two torsion free connections $\nabla^1$ and $\nabla^2$ are projectively equivalent if there exist a 1-form $\beta$ such that
\begin{eqnarray}
\nabla^1_XY-\nabla^2_XY=\beta(X)Y+\beta(Y)X.
\end{eqnarray}
A class of projectively equivalent connections defines a {\em projective structure} on $M$. 
Projective structures can be induced by local diffeomorphisms. 

\subsubsection{Projectively flat structures}
A projective structure is {\em flat} if every point has a neighbourhood on which the projective structure is given by a torsion free flat connection.

The projective space   $\rpn$ admits a projectively flat structure given  by the  affine charts. Conversely, a manifold $M$ of dimension $n$ is equipped with a \sp, if there exist
\begin{itemize} 
\item a representation $\rho$ --the {\em holonomy representation}-- of $\pi_1(M)$, the fundamental group of $M$, with values in the projective group  $PSL(n+1,\RR)$~; 
\item a local diffeomorphism $f$, the {\em developing map}, of the universal cover $\wtm$ of $M$ with values in $\rpn$, which is  $\rho$-equivariant,  that is which satisfies
$$
\forall\, x\in\wtm, \, \forall\, \g\in\pi_{1}(M),\quad f(\g x)=\rho(\g)f(x).
$$
\end{itemize}
The structure    \ on $M$ is the one induced by the projective structure on  $\rpn$ by $(f,\rho)$.

In other words, a  \sp\ on a manifold  is nothing else that a 
$(\rpn, PSL(n+1,\RR))$-{\em structure }, or, in short, a $\rpn$-{\em structure }. We shall restrict ourselves in this article to structures such that the holonomy representation lifts to $SL(n,\mathbb R)$.

\subsubsection{Convex structures}

A $\rpn$-structure  is  {\em convex}  if the developing map is a homeomorphism to a convex set in  $\rpn$. It is {\em properly convex } if this  convex  set is included in a compact convex set of an affine chart.

\subsection{Projective structures and connections}\label{sec:projconn}
We explain and relate  in this section two different points of view on projectively flat structures.\begin{itemize}
\item In the first paragraph, we explain that pairs consisting of a flat connection and a section of a rank $n+1$ vector bundle over an $n$-dimensional manifold can give  rise to flat projective structures.
\item In the second paragraph, we explain that the pair consisting of a torsion free connection and a symmetric tensor - satisfying some compatibility condition - on a  manifold also gives to a flat structure   
\end{itemize}
\subsubsection{Flat connections and sections}

We consider the trivial bundle  $E_M=M\times E$ where $E$ is an $n+1$-dimensional vector space equipped with a volume form $\omega$. Let $\nabla$ be a connection on $E_M$ preserving the volume form. We observe that each section $u$ of $E_M$, defines
an element of $\Lambda^n(TM^*)$, given by
$$
\Omega^0_u(X_1,\ldots,X_n)=\omega(\nabla_{X_1}u,\ldots,\nabla_{X_n}u, u).
$$
We say $u$ is {\em $\nabla$-immersed} if i$\Omega^0_u$ is non degenerate.

We now relate these notions to projective structures. 
If $\nabla$ is flat, we identify  $E_M$ on the universal cover $\tilde M$ on $M$ with $E\times \tilde M$ so that $\nabla$ is the trivial connection. Let $\rho$ be the holonomy representation of $\nabla$. A section $u$ of $E_M$ is then identified with a $\rho$-equivariant map from $\tilde M$ to $E$.  We denote by $x\to [x]$ the projection from $E\setminus\{0\}$ to $\mathbb P (E)$.
A section $u$ is $\nabla$-immersed, if it is a non zero section and if the associated  $\rho$-equivariant map $[u]$ from $M$ to $\mathbb P (E)$ is an immersion.

It follows that a pair $(\nabla,u)$ such that $\nabla$ is flat and $u$ is $\nabla$-immersed gives rise to a flat projective structure. Conversely it is immediate to check every flat projective structure whose holonomy lifts to $SL(n+1,\mathbb R)$ is obtained this way, maybe after going to a double cover.

We obtain immediately from the previous construction the following proposition
\begin{proposition}\label{fu}
If $(\nabla_1,u_1)$ and $(\nabla_2,u_2)$ give rise to two flat projective structures equivalent by a diffeomorphism $\phi$, then there exist
\begin{itemize}
\item a bundle automorphism $\Phi$ of $E_M$ over $\phi$,
\item a non zero function $f$,
\end{itemize}
such that $\Phi^*\nabla_1=\nabla_2$ and $\Phi^*u_1=f u_2$.
\end{proposition}

\subsubsection{Connections and symmetric tensors }\label{connsym}

Let  $\nabla$ be a torsion free  connection on  $M$. Let $h$ be a symmetric two-tensor on $M$. Let  $L$ be the trivial bundle  $\RR\times M$. We associate to the pair $(\nabla,h)$ a connection  $\nabla^h$ on $TM\oplus L$ given by 
$$
\nabla^{h}_{X}\vecteur{Z\\ 
\la} = \vecteur{\nabla_{X}Z+\la X\\ 
L_X\la+h(Z,X)}
$$
We say $(\nabla, h)$ satisfy {\em Condition (E)} if ,
\begin{itemize}
\item $\nabla$ preserves a volume form,
\item $\nabla^h$ is flat.
\end{itemize}
Note that if $\nabla$ satisfies Condition (E), then $\nabla^h$ preserves a volume form on $TM\oplus L$. Moreover, the conditions that $\nabla$ is torsion free and $h$ symmetric are redundant. Finally, $\nabla^h$ is flat is equivalent to the following two equations 
\begin{eqnarray}
\left\{\begin{array}{rcl}
d^\nabla h (X,Y,Z)=\nabla_X h(Y,Z)-\nabla_Y h(X,Z)Y&=&0\\ \label{condE}
R^\nabla(X,Y)Z-h(X,Z)Y+h(Y,Z)X&=&0.
\end{array}
\right.
\end{eqnarray}
We now relate this to the previous paragraph. 
Let as above $\bar\nabla$ be a connection on $E_M$ preserving a volume form $\omega$. We observe that each section $u$ of $E_M$, defines
a  tensor element $\Omega^2_u$ of $\Lambda^n(TM^*)\otimes S^2(TM)$, by the following formula
$$
\Omega^2_u(X_1,\ldots,X_n, Y,Z)=\omega(\bar\nabla_{X_1}u,\ldots,\bar\nabla_{X_n}u,\bar\nabla_{X}\bar\nabla_{Y} u).
$$
If $u$ is $\bar\nabla$-immersed, then we can write
$$
\Omega^2_u=\Omega_u^0\otimes S_u.
$$
Note  that the symmetric two tensor $S_u=S^{\bar\nabla}_u$ is independent of the choice of $\omega$.

We observe that if $u_0$ is the canonical section of $TM\oplus L$, then $h=S^{\nabla^h}_{u_0}$. Conversely if $(\nabla^0,u)$ is a pair such that $u$ is $\nabla$-immersed, then the following mapping is an isomorphism
\begin{eqnarray*}
\Phi:\mapping{TM\oplus L}{E_M}{(X,\lambda)}{\nabla_X u+\lambda u}.
\end{eqnarray*}
Moreover there exists a connection $\nabla=\nabla^u$ on $M$ such that 
$$
\Phi^*\nabla^0=\nabla^h, \hbox{ with } h=S^{\nabla^0}_u.
$$
 
This construction is related to projective structures by the following result.

\begin{proposition}
If $(\nabla,h)$ satisfies Condition (E), then $\nabla$ is projectively flat and the holonomy of the the corresponding structure is the holonomy representation of $\nabla^h$. Conversely, every flat projective structure on $M$, whose holonomy lifts to $SL(n+1,\mathbb R)$ is obtained this way. 
\end{proposition}

\proof the only point to be proved is that $\nabla$ is projectively flat and defines the same projective structure. Let $\gamma(t)$ be  a geodesic for $\nabla$, then the sub-bundle 
$$
P=\mathbb R\dot\gamma(t)\oplus \mathbb R \subset TM\oplus L\mathbb R,
$$ 
is parallel along $\gamma(t)$. Therefore, using the local trivialisation given by $\nabla^h$, $[u(\gamma(t)]$ is the projective line defined by $P$. \qed 
\vskip 0.5truecm

We also prove an independent proposition which will proved technically useful in the sequel.

\begin{proposition}\label{suD}
We consider the connexion $\nabla^u$ on $TM$, such that 
$$
\nabla_X\nabla_Y u= S_u(X,Y)\cdot u + \nabla_{\nabla^u_X Y}u.
$$
Then

\begin{eqnarray}
S_{\frac{u}{f}}&=&S_u -\frac{(\nabla^u)^2f}{f},\label{sfusu}\\
\nabla^{\frac{u}{f}}_XY&=&\nabla^u_X Y-\frac{df(X)}{f}Y-\frac{df(Y)}{f}X\label{nfunu}.
\end{eqnarray}

\end{proposition}
\proof We  consider $\pi_u : E_M\rightarrow TM$, such that 
$$
\pi_u(u)=1, \ \ \pi_u(\nabla_X u)=0.
$$
Then 
$$
0=\pi_{\frac{u}{f}}(\nabla_X (\frac{u}{f}))=\frac{1}{f}(\pi_{\frac{u}{f}}(\nabla_X u)-df(X)).
$$
Thus,
$$
\pi_{\frac{u}{f}}(\nabla_X u)=df(X).
$$
Also
\begin{eqnarray*}
\nabla_X\nabla_Y\frac{u}{f}&=&\frac{1}{f}\nabla_X\nabla_Y u -\frac{1}{f^2} L_XL_Y f u\\ &+&2\frac{df(X)df(Y)}{f^3}u-\frac{df(X)}{f^2}\nabla_Y u --\frac{df(Y)}{f^2}\nabla_X u
\end{eqnarray*}
By definition, $S_u(X,Y)=\pi_u(\nabla_X\nabla_Y u)$.
Then
\begin{eqnarray*}
S_{\frac{u}{f}}(X,Y)&=&\pi_{\frac{u}{f}}(\nabla_X\nabla_Y \frac{u}{f})\\
&=&-\frac{1}{f} L_XL_Yf+\frac{1}{f}\pi_{\frac{u}{f}}(\nabla_X{\nabla_Y}u)\cdot u)\\
&=&-\frac{1}{f} L_XL_Yf+\frac{1}{f}\pi_{\frac{u}{f}}(\nabla_{\nabla^u_X Y}u)+S_u(X,Y)\\
&=&-\frac{1}{f}(L_XL_Yf-df(\nabla^u_X Y))+S_u(X,Y)\\
&=&-\frac{1}{f}(\nabla^u)^2_{X,Y}f +S_u(X,Y).
\end{eqnarray*}
Finally, writing $v=\frac{u}{f}$, we have
\begin{eqnarray*}
\nabla_{\nabla^{v}_XY}\frac{u}{f}&=&\nabla_X\nabla_Y \frac{u}{f}-S_{v}(X,Y)\frac{u}{f}\\
&=&-\frac{df(X)}{f}\nabla_Y\frac{u}{f}-\frac{df(Y)}{f}\nabla_X\frac{u}{f}+\nabla_{\nabla^u_X Y}\frac{u}{f}\\
\end{eqnarray*}
It follows that
$$
\nabla^{\frac{u}{f}}_XY=\nabla^u_X Y-\frac{df(X)}{f}Y-\frac{df(Y)}{f}X.
$$
\qed

\section{Convex projective structures}

In this section, we explain how convex projective structure can be interpreted using the point of view of the previous Section.  
Our main result is Theorem \ref{rpncon}.

\subsection{Hypersurfaces  and convex  $\rpn$-structures}\label{rpnS}

Let  $E$ be a vector space of dimension $n$. We say a an immersed hypersurface $S$ in $E\setminus\{0\}$ 
\begin{itemize}
\item is {\em locally convex}, if  every point in $S$ has a neighbourhood $U$ in $S$ which is a subset of the boundary of a convex set. 
\item  is {\em locally strictly convex},  if furthermore  $U$ does not contain any segment. In particular a locally strictly convex separates $E$ locally in two connected components, the {\em interior} --which coincides with the convex set -- and an {\em exterior}
\item We say a  strictly locally convex hypersurface is {\em radial}, if the radial vector -- the vector pointing from the origin -- points inward.
\item Finally, we observe that if  $S$ is  strictly locally convex, radial and properly immersed, it bounds a convex set which does not contain the origin.
\end{itemize}

Every such hypersurface admits a natural properly convex  structure. Indeed the projection from 
 $S$ to $\PP(E)$ is an immersion. Since $S$ is strictly convex, 
its projection is  a convex  set whose closure is compact and  included in an affine chart.

Conversely, the following lemma, basically due to Vinberg,  shows every properly convex  structure is obtained this way. 

\begin{lemma}\label{convS}
{ Let  $M$ be  manifold equipped with a properly convex  projective structure given by the pair $(f,\rho)$. Then there exists a proper $\rho$-equivariant immersion  $g$ of $\wtm$ whose image is strictly convex and radial and and such that $\pi\circ g = f$ where $\pi$ is the projection 
$E\setminus\{0\}$ on  $\PP(E)$.}
\end{lemma}
  
\proof Let $\cC$ be an open  proper convex set of $\PP(E)$. Let $C^0$ be the  cone of $E$ obtained from $\cC$
$$
C^0=\pi^{-1}(\cC).
$$
Let $C$ be the convex cone which is one the two connected component of $C^0$. Let $C^*$be the dual cone of $C$. Let $df$ be a volume form on the dual $E^{*}$ of $E$. Let $V$ be the characteristic function of 
 Vinberg [Vi]
defined by
$$
\forall x\in C, \ V(x) = \int_{C^{*}}e^{-f(x)}df.
$$
This function $V$ is convex and one shows easily that the   
hypersurface $S=V^{-1}\{1\}$ is properly immersed strictly convex and radial and $\pi$ is a diffeomorphism  from $S$ to $\cC$.

Since  $V$ is invariant by the subgroup  of the special affine group that leaves 
$C$ invariant, we obtain that  $S$ is also invariant.

In particular, if  $f$ is the developing map of a properly convex structure. We can define
$$
g:
\mapping{\wtm}{S}{x}{\pi^{-1}( f(x)).}
$$
The mapping $g$ satisfies the condition of the lemma. \qed

\subsection{Convex $\rpn$-structures and connections}

We explain in this section that a properly convex  $\rpn$-structure 
on a compact manifold $M$  is equivalent to the data  of a connection $\nabla$ on $TM$, and a positive definite metric $g$ on $TM$ such that $(\nabla,g)$ satisfies Condition (E).  We use the language of Paragraph \ref{connsym}

\begin{theorem}\label{rpncon}
If $(\nabla,g)$ satisfies Condition (E) and $g$ is positive definite,  then $\nabla$ is projectively  flat and defines a properly convex structure on $M$. Conversely, every properly convex projectively flat structure on $M$ is obtained this way. 
\end{theorem}

We first observe that since $\nabla^{g}$ is flat, the bundle $T\wtm \oplus L$ over $\wtm$ is isomorphic to the trivial bundle $E\times \wtm$. Let $p$ be the projection of $T\wtm \oplus L$ to $E$.
Let $\rho$ be the holonomy representation of the flat connection $\nabla^g$. Let $s: m\rightarrow (0,1)$ be the canonical section of $L\subset T\tilde M\oplus L$. Let $\phi=p\circ s$. By construction $\phi$ is a $\rho$-equivariant mapping of $\wtm$ in $E$.

We  prove the Theorem in several steps
\begin{proposition}\label{fimmer}
The mapping $\phi$ is an immersion. 
\end{proposition}
\proof
We have by construction
$$
D\phi (X)=p(\nabla^g_X s)=p(X,0).
$$
It follows $\phi$ is an immersion. \qed
\begin{proposition}\label{fconv}
The immersed  hypersurface $\phi(\wtm)$ is strictly locally convex and radial.
\end{proposition}
\proof Let $S$ be a hypersurface in $E$. For every $s\in S$, let 
$\pi_s$ be the projection $E\rightarrow E/T_s S$. Let $D$ be the flat connection in $E$. For every vector field $X$ on $S$, we recall that there exists asymmetric 2-tensor $B$ such that
$$
\pi_s(D_X,X)=B(X_s,X_s).
$$
A  hypersurface is strictly convex at $s$ if and only if  for every non zero  vector  $X$  in $T_sS$, we have
$$
B(X_s,X_s)\not=0.
$$
In our case, let $S=\phi(\wtm)$, then 
$$
B(X_s,X_s)=g(X_s,X_s)>0.
$$
This also prove that $S$ is radial.
\qed
\vskip 0.5truecm
The next proposition is less straightforward.
\begin{proposition}\label{fproper}
$\phi$ is proper.
\end{proposition}
\proof We first note that  the geodesics for $\nabla$ are precisely those curves whose image by $\phi$ lies on 2-plane which passes through the origin. 
Let  $x_{0}$ a  point in $\wtm$ and  $U$ the domain of  $T_{x_{0}}\wtm$ on which the exponential map $\exp$ for $\nabla$ is defined. Let $v$ be a non zero vector in $T_{x_{0}}\wtm$. Let $I=]a,b[$ be the maximal interval for which 
$Iv\subset U$.   Let $c_0(t)=\exp(tv)$.  We now prove that
$$
\psi_v
\mapping{I}{E}{t}{\phi(\exp(tv))}
$$
is proper. We know that $\psi_v(I\cap U)$ is a strictly convex planar curve $c$ which is furthermore radial. To prove that 
$\psi_v$ is proper it suffices to show that the length of this curve is infinite for a euclidean metric $\langle\ ,\rangle$ on $E$. 
Let $R:x\rightarrow x$ be the radial vector field of $E$. Let $\mu (t)=\sqrt{g(\dot{c_0}(t),\dot{c_0}(t))}$.
We remark that 
\begin{eqnarray*}
\ddot{c}&=&\mu^2 R, \\ \dt{R}&=&\dot{c}.
\end{eqnarray*}

Let $\la (t) =\sqrt{\langle \dot{c}(t),\dot{c}(t)\rangle}$.
We have
\begin{eqnarray}
\dt{\la}&=&\frac{\langle \dot{c},\ddot{c}\rangle}{\lambda}\cr
&=&\mu^2.\frac{\langle \dot{c},R\rangle}{\la}
\end{eqnarray}
However, $t \rightarrow \frac{\langle \dot{c}(t),R(t)\rangle}{\la(t)}$
 is increasing:
\begin{eqnarray*}
\dt{\ }\frac
{\langle \dot{c},R \rangle }
{\la}
&=&
\frac{
\la(\langle \ddot{c},R \rangle +\langle \dot{c},\dot{c}\rangle).
-\mu^2\la^{-1}\langle R,\dot{c}\rangle^2
}{\la^2}
\\
&=&\frac{\mu^2\la \langle R,R \rangle +\la\langle \dot{c},\dot{c}\rangle
-\mu^2\la^{-1}\langle R,\dot{c}\rangle^2
}{\la^2}
\\ 
&=&
\frac{(\langle R,R \rangle \langle \dot{c},\dot{c}\rangle -\langle R,\dot{c}\rangle^2)\mu^2\la^{-1}}{\la^2}
+\frac{\langle \dot{c},\dot{c}\rangle}{\la}
\\
&\geq& 0.
\end{eqnarray*}
Let us choose an euclidean metric $g_{1}$ such that
$\langle \dot{c}(0),R \rangle  > 0$.
Thus
$$
\forall  t > 0, \frac{\langle \dot{c}(t),R \rangle }{\lambda (t)}\geq K_{1}>0.
$$
Hence,
$$
\forall t > 0, \ \  \dt{\la}\geq K_{1}\mu^2.
$$
But, by Lemma \ref{complete} of the Appendix.
$$
\forall t , \  \vert\dt{\mu}\vert\leq  K_{2}\mu^{2} .
$$
Thus 
$$
\forall t>  0, \ \  \dt{\la} \geq K_{3}\dt{\mu}.
$$
It follows again by Lemma \ref{complete} of the Appendix that
$$
\int_{0}^{b}\lambda(s)ds\geq K_3\int_{0}^{b}\mu(s)ds +K_4=+\infty.
$$
Finally, we show that $\phi\circ\exp$ is proper. Let $\{u_n\}\inn$ be a sequence of points in $U$ such that $\{\f\circ\exp(u_{n})\}\inn$ converge to $z_{0}$. By choosing a subsequence, we may suppose that the sequence of lines $\{D_{n}=\RR u_{n}\}\inn$
converges to $D_{0}$, and that - according to the previous discussion -  $\{\f\circ\exp(D_{n}\cap
U)\}\inn$ converges  to a locally convex  \ radial curve  $c$. But then, 
$\phi\circ\exp(D_0\cap U)$ is an open set in $c$, and since it has infinite length it coincides with  $c$. Thus,  $\{u_{n}\}\inn$ converges to a point $u_{0}$ in  $D_0\cap U$ such that  $\phi\circ\exp(u_{0})=z_{0}$. 
\qed
\subsubsection{Proof of Theorem \ref{rpncon}}
We prove the first part of the Theorem. By Proposition \ref{fconv}, \ref{fimmer} and \ref{fproper}, $\phi$ is a proper immersion and $\phi(M)$ is a locally convex proper hypersurface. By Section \ref{rpnS}, we can construct a flat projective structure on $M$, whose projective lines are the intersection of 2-planes with $\phi(M)$. But the latter are geodesics on $\wtm$ for $\nabla$, therefore $\nabla$ is projectively flat.

Conversely, by Lemma \ref{convS}, any properly convex projectively flat  structure on $M$ can be induced by  proper $\rho$-equivariant immersion in $E$ of $\tilde M$ whose image is  a locally strictly convex hypersurface $S$. Since $S$ is radial, we have the decomposition
$$
TE\vert_S=TS\oplus \mathbb R R.
$$
The flat connexion $\nabla^0$ on $TE\vert_S$ decomposes therefore, if $Z$ is a vector field on $S$ as
$$
\nabla^{0}_{X}\vecteur{Z+ 
\la R} = \nabla_{X}Z+\la \nabla_X R+ 
L_X\la.R+g(Z,X).R.
$$
We observe that  $\nabla^0_X R=X$, and $g(X,X)>0$ since $S$ is strictly locally convex and radial. It follows that $(\nabla,g)$ satisfies Condition  (E) \ref{condE}.

\section{Cubic  holomorphic differentials  and convex  $\rpd$- structures.}

The two main results of this section, Theorem \ref{cub->rp} and \ref{rp->cub}, provide a bijection of the space of pairs consisting of a complex structure and  a cubic holomorphic differential, with the moduli space of convex projective structures.

\subsection{From cubic  holomorphic differentials  to convex  $\rpd$- structures.}

The content of this section is another version of Wang Chang Pin result \cite{wcp} which  states that every cubic holomorphic differential on a compact surface gives rise to an affine sphere. The proofs is slightly different. However for the sake of completeness, and to make this article as much self contained as we can,  we recall the construction. 
\begin{theorem}\label{cub->rp}
Let $\omega$ be a cubic holomorphic differential. For any metric $g$, let $A_g$ be the element of 
of  $TS^{*}\otimes {\rm End}( TS)$ defined by  
$$
\Re(\omega(X,Y,Z))=g(A_g(X)Y,Z).
$$
Then there exists a unique metric $g$ in the conformal class of $J$ such that if $\bar\nabla$ is the Levi-Civita connexion of $g$, then
$(\nabla=\bar\nabla +A_g,g)$ satisfies Condition  $(E)$. Moreover the area for $g$ is parallel for $\nabla$. \end{theorem}

Let $g$ be a metric conformal to $J$. Let   $\nabla$ be the Levi-Civita connection of $g$. Let $A\in\Omega^1(S)\otimes TS^*\otimes TS$.  Then a straightforward check shows that $(\nabla +A,g)$
satisfies Condition $(E)$ if and only if
\begin{enumerate}
\item $A(X)$ is symmetric and trace free
\item $A(X)Y = A(Y)X$
\item  $d^{\nabla} A=0$
\item \label{eqcourb}$
R^{g}(X,Y)Z + [A(X),A(Y)]Z + g(Y,Z)X-g(X,Z)Y=0.
$
\end{enumerate}
Let now $A=A_g$ be associated to $g$ and $\omega$  as in the statement of the theorem. The first three conditions are satisfied since $\omega$ is a cubic holomorphic differential. We concentrate on the last condition. 
We first rewrite it in a more compact way. There exist a positive definite quadratic form   $\cG^{g}$, depending on $g$ and defined on the space of cubic differentials such that 
$$
[A(X),A(Y)]Z=-\cG^{g}(\omega,\omega)(g(Y,Z)X-g(X,Z)Y).
$$
Note  that
$$
\cG^{\la g}=\frac{1}{\la^3}\cG.
$$
We can rewrite (\ref{eqcourb}) as  
\begin{eqnarray}
k_{g}+1- \cG^{g}(\omega,\omega)=0\,  \label{eqcourb2}
\end{eqnarray}
where $k_g$ is the curvature of $g$. 
Write $g=\la g_{0}$, where  $g_{0}$ has constant curvature -1. Let $f=\cG^{g_{0}}(\omega,\omega)$ and  $\mu=\frac{1}{2}\log (\la)$. Recall that
$$
k_{g}=e^{-2\mu}(\D\mu -1).
$$
Then  we rewrite  (\ref{eqcourb2}) as
\begin{eqnarray}
f - e^{4\mu}\D\mu-e^{6\mu} + e^{4\mu} = 0.
\end{eqnarray}
Therefore the result follows from
\begin{lemma}\label{lemlaplace}
Let $S$ be a compact hyperbolic surface. Let $f$ be a positive function on $S$. Then there exists a unique function $\mu$ such that
$$
e^{4\mu}\D\mu+e^{6\mu} - e^{4\mu} = f
$$
where $\D=\frac{\pr}{\pr X_{1}^{2}}+ \frac{\pr}{\pr X_{2}^{2}}$ is the Laplacian
\end{lemma}
The rest of this paragraph is devoted to the proof of this Lemma.
Let 
$$
H(\mu)=e^{+4\mu}\D\mu+e^{6\mu}-e^{4\mu}
$$
\subsubsection{A-priori estimates}
Let $\mu$ be a $\ci$-function on $S$. Let $f=H(\mu)$. We assume that $f\geq 0$. We want to control $\mu$ using $f$. Our main result is the following
\begin{lemma}\label{apriori}
For every $A$, there exists $B$ such that $\Vert f\Vert_{C^1}\leq B$ implies
$\Vert \mu \Vert_{C^1}\leq A$
\end{lemma}
\proof
We first obtain $C^0$ estimates using the maximum principle.
At a minimum of $\mu$, we have $\Delta\mu \geq 0$. Thus
$
f\leq e^{6\mu}-e^{4\mu}
$.
Therefore $e^{6\mu}-e^{4\mu}\geq 0$. Hence $\mu \geq 0$.
Now at  a maximum of $\mu$,   we have $\Delta\mu \leq 0$ thus
$
f\geq e^{6\mu}-e^{4\mu}
$.
This proves  that
for every $A$, there exists $B$ such that $\Vert f\Vert_{C^0}\leq B$ implies
$\Vert \mu \Vert_{C^0}\leq A$.

Next we prove the $C^1$ estimates. Since $g_0$ has constant curvature -1, we have
$$
\D\|d\mu\|^{2}= 2\langle d(\D\mu),d\mu \rangle
+\|\nabla^{2}\mu\|^{2}-2\|d\mu\|^{2}\, .
$$
Let $f=H(\mu)$, thus:
\begin{eqnarray*}
\D\|d\mu\|^{2}&=&\|d\mu\|^{2}(8e^{-4\mu}f+4e^{2\mu}-2)+\|d^{2}\mu\|^{2}-2e^{-4\mu}\langle
 df, d\mu\rangle
\\
&\geq&\big(\|d\mu\|^{2}(4f +2e^{6\mu}-e^{4\mu})-\|df\|\cdot\|d\mu\|\big)2e^{-4\mu}.
\end{eqnarray*}
We have shown that $\mu\leq 0$. Thus
$$
4f+2e^{6\mu}-e^{4\mu}>1.
$$ 
Therefore when $\Vert d\mu\Vert$ is maximum then
$$
0\geq \D\|d\mu\|^{2}\geq \big(\|d\mu\|^{2}-\|df\|\cdot\|d\mu\|\big)2e^{-4\mu}.$$
Thus
$$
\Vert d\mu \Vert \leq  \sup \Vert df \Vert.
$$
\qed
\subsubsection{Proof of Lemma \ref{lemlaplace}}
Let 
$$\cF=\{\mu/H(\mu) \geq 0\}$$
We first want to prove that 
$$
H: \cF\to \CI(S,[0,+\infty[)
$$
is a homeomorphism. It suffices to show that
\begin{enumerate}
\item $H$ is a local homeomorphism,\label{Hi}
\item $H$ is proper,\label{Hii}
\item $\cF$ is connected.\label{Hiii}
\end{enumerate}

{\em We prove (\ref{Hi}).} Since $H$ is an elliptic operator, it suffices to show by the local inversion theorem that the linearised operator  $L^{H}_{\mu}$ of $H$ at  $\mu$ is invertible.
A straightforward computation yields 
$$
L^{H}_{\mu}(\la)=4\la H(\mu)+2\la e^{6\mu}-e^{4\mu}\D\la\, .
$$
By hypothesis
$2H(\mu)+e^{6\mu}>0$.  By the maximum principle we deduce that  the kernel $L^{H}_{\mu}$ is reduced to zero. Since the index of $L^H_\mu$ is zero, it follows  $L^H_\mu$ is invertible

{\em We prove (\ref{Hii}).} Assume that $\{H(\mu_n)\}\inn$ converges. It follows from Lemma \ref{apriori}, that $\{\Vert\mu_n\Vert_{C^1}\}\inn$ is bounded. Then classical arguments shows that $\{\mu_n\}\inn$ converges. For the sake of completeness, we give a slightly non standard proof of that fact  in Appendix \ref{lapla}.

{\em We prove (\ref{Hiii}).} For the moment we have proved that $H$ is a finite index covering. To prove that $\cF$ is connected, it suffices to show that for some $f$, there is a unique solution of $H(\mu)=f$. Indeed, for $f=0$, the equation $H(\mu)=0$ says that  $e^{2\mu}g_{0}$ has constant curvature -1. Hence $\mu=0$. \qed

\subsection{From convex $\rpd$-structures  to cubic  holomorphic differentials}

We prove the following theorem

\begin{theorem}\label{rp->cub}
Let $S$ be  equipped with a  convex  projective structure. Then there exist a unique pair
$(\nabla,g)$ defined up to diffeomorphism, such that 
\begin{itemize}
\item $\nabla$ defines the projective structure,
\item $g$ is definite positive,  
\item if $\omega_g$ is the volume form of $g$, then $\nabla \omega_g=0$.
\end{itemize}
Moreover, in this situation, let $\bar\nabla$ be the Levi-Civita connexion of $g$, then
$$
\Omega(X,Y,Z)=g(\nabla_X Y-\bar\nabla_X Y,Z).
$$
is the real part of a cubic holomorphic differential.
\end{theorem}

We can translate the theorem in the the  language of splitting using the notations of Paragraph \ref{connsym}.

\begin{theorem}{\sc[Splitting] }
\label{theo:split}
Let $S$ be  equipped with a  convex  projective structure. Let $E$ be the associated flat vector bundle equipped with a parallel volume form $\Omega$, then there exist a unique splitting $E=P\oplus\mathbb Ru$ such that
\begin{itemize}
\item $\nabla^u$ defines the convex structure
\item $\forall X,\ \  \nabla_X u\in P$.
\item the volume form of the quadratic form $S^u$ is $\Omega(\nabla_X u,\nabla_Y u, u)$
\end{itemize}
\end{theorem}

\subsubsection{Back to cubic holomorphic differential}
We first prove the second part of the Theorem.
\begin{proposition} Let  $(\nabla, g)$ be such that 
\begin{itemize}
\item $\nabla^g$ is flat
\item
$\nabla$ preserves the area form of $g$.
\end{itemize}
Let $\nabla_g$ is the Levi-Civita connexion of $g$.
Then
$$
\Omega(X,Y,Z)=g(\nabla_{X}Y-(\bar\nabla)_{X}Y,Z)
$$
is the real part of a cubic holomorphic differential.
\end{proposition}
  
\proof Let  $A(X)Y = \nabla_{X}Y-\bar\nabla_{X}Y$. A straightforward check shows that $\Omega$ is the real part of a cubic holomorphic differential if and only if 
\begin{enumerate}
\item $d^{\bar\nabla} A=0$\label{dnabla}
\item  $A(X)Y=A(Y)X$\label{torsion}
\item  $A(X)$ is symmetric.\label{symme}
\item $A(X)$ is trace free.\label{trace}
\end{enumerate}
Since $\bar\nabla$ preserves the area form of  $g$, $A(X)$ is trace free . Thus Condition (\ref{trace}) is satisfied. Recall that by definition $\nabla^g$ is given by
$$
\nabla_{X}^{g}\vecteur{Z\\  \la} = \vecteur{\nabla_{X}Z +\la X\\  L_{X}\la 
+g(Z,X)}.
$$
We compute various parts of the curvature tensor $R$ of $\nabla^g$. First, 
$$
0=R(X,Y)\vecteur{0\\  1},$$
yields
$$
\nabla_X Y-\nabla_Y X-[X,Y]=0.$$
 Thus $\nabla$ is torsion free. Hence, Condition (\ref{torsion}) is satisfied.
Let us compute the other part of the curvature tensor 
$$
0=R(X,Y)\vecteur{Z\\  0},$$
yields
$$
\left\{
\begin{array}{rcl}
R^{\nabla}(X,Y)Z + g(Y,Z)X -
g(X,Z)Y&=&0\\ 
g(A(Y)Z,X)-g(A(X)Z,Y)&=&0
\end{array}
\right .
$$
Using Condition (\ref{torsion}),  the second line of this equation reads
$$
g(A(Z)Y,X)=g(A(Z)X,Y).
$$ Thus, Condition (\ref{symme}) is satisfied. 
The first line yields
$$
R^{g}(X,Y)Z+d^{\nabla^g} A(X,Y)Z + [A(X),A(Y)]Z + g(Y,Z)X-g(Z,X)Y=0.
$$
However, the linear operator $Z\mapsto g(Y,Z)X-g(X,Z)Y$ is antisymmetric. The same holds for  $R^g(X,Y)$ and $[A(X),A(Y)]$, since $A(X)$ is symmetric. Thus
 $Z\mapsto d^{\nabla}A(X,Y)(Z)$ is antisymmetric. However,  $A(X)$ is symmetric, thus
  $Z\mapsto d^{\nabla}A(X,Y)(Z)$ is symmetric. Therefore
  $d^{\nabla}A=0$ and Condition (\ref{dnabla}) is satisfied.\qed

\subsubsection{Proof of Theorem \ref{rp->cub}}

Let $\nabla^0$ be a flat connection on $E_M=M\otimes \mathbb R^{n+1}$ preserving a volume form $\Omega^0$. Let  $u$ be a section of $E_M$ such that $(\nabla^0,u)$ give rise to the projective structure. It follows from Paragraph 
\ref{connsym} that for every pair $(\nabla,g)$ where 
\begin{itemize}
\item
$g$ is positive,
\item $\nabla$ equivalent to the convex structure, 
\item $\nabla$ preserves a  volume form,
\item  $\nabla^g$ is flat, 
\end{itemize}
there exists a  function $f$ such that 
$$
S_{\frac{u}{f}}=g.
$$
Let now 
$$
\mathcal U=\{f\in C^\infty(S), \hbox{  such that } S_{ \frac{u}{f}} \hbox{ is positive definite }\}.
$$
For every $f\in\mathcal U$,
\begin{itemize}
 \item Let $\omega^f$ be the area form of $S_{\frac{u}{f}}$. 
 \item Let $\nu^f(X,Y)=\Omega^0(\nabla_X \frac{u}{f}, \nabla_Y \frac{u}{f}, \frac{u}{f})$
 \end{itemize}
The theorem will follow from the following proposition.

\begin{proposition}\label{nu}
The operator 
$$
D:\mapping{\mathcal U}{C^\infty(S,]0,\infty[)}{f}{\frac{\nu^f}{\omega^f}},
$$
is a diffeomorphism.
\end{proposition}

\subsubsection{Proof of Theorem \ref{rp->cub} from Proposition \ref{nu}.}
Let $\nabla=\nabla^u$ be the connection on $S$ given by
$$
\nabla^0_X\nabla^0_Y u=S_u(X,Y)u+ \nabla^0_{\nabla_XY}u.
$$
We first recall from Equation (\ref{sfusu}) of Proposition \ref{suD}, that
$$
S_{ \frac{u}{f}}=S_u-\frac{\nabla^2 f}{f}.
$$
In particular, $S_{\frac{u}{kf}}=S_{ \frac{u}{f}}$.
Hence $D(kf)=k^{3}D(f)$. Now observe that $\omega_f$ is parallel for $\nabla^{ \frac{u}{f}}$. It follows that the pairs $(\nabla,g)$ for which $\nabla\omega_g=0$, corresponds to the functions $f$, up to a multiplicative constant, such that $D(f)$ is constant. Thus, the first part of Theorem \ref{rp->cub} follows from Proposition \ref{nu}.

\subsubsection{Proof of Proposition \ref{nu}}
 
Let $A(f)$ be the symmetric endomorphism defined by  
$$
\nabla^{2}f(X,Y)=S_u(A(f)X,Y).
$$Observe that $\nabla=\nabla_u$ is torsion free.
Thus $A(f)$ is symmetric. Then
\begin{eqnarray}
\nu^{f} &=&\det\big(1-\frac{A(f)}{f}\big)\nu^{1}\\ 
&=&\det\big(1-\frac{A(f)}{f}\big)f^{3}D(1)\omega^{f}
\end{eqnarray}
Therefore
$$
D(f)=\det\big(f^{3/2} - f^{1/2}A(f)\big)D(1)\, .
$$
The proof of the Proposition will follows from the following three steps
\begin{enumerate}
\item $D$ is a local diffeomorphism from $\cU$ in $\CI(S,]0,+\infty[)$,
\item $D$ is proper,
\item $\cU$ is connected.
\end{enumerate}
This will require three separate propositions.

\begin{proposition}
$D$ is a local diffeomorphism from $\cU$ in $\CI(S,]0,+\infty[)$,
\end{proposition}
  
\proof We first compute the linearised operator $L^{D}_{f}$ at $f$. Let 
$$G(f)=f^{3/2}-f^{1/2}A(f).$$
Then a straightforward computation give
\begin{eqnarray*}
L^{D}_{f}(\mu) &=&\frac{D(1)}{D(f)}\tr\big((\frac{3}{2}f^{\frac{1}{2}} \mu -\frac{1}{2}
f^{-\frac{1}{2}}A(f)\mu-f^{\frac{1}{2}}A(\mu))\circ G(f)^{-1}\big)\\ 
&=&\frac{D(1)}{D(f)}\tr\big(\mu (f^{\frac{1}{2}}+
\frac{G(f)}{f})-f^{\frac{1}{2}}A(\mu))\circ G(f)^{-1}\big)
\end{eqnarray*}

Since $f\in\cU$, $G(f)$ is a positive symmetric operator, it follows that $L^{D}_{f}$ is an elliptic operator. By the implicit function theorem, to prove the Proposition it suffices to show that $L^{D}_{f}$ is invertible. Since
$L^{D}_{f}$ is homotopic to a Laplacian, its index is zero. It thus suffices to show that 
$L^{D}_{f}$ is injective.
Let $\mu$ so that  $L^{D}_{f}(\mu)=0$. Hence 
\begin{eqnarray}
\tr(A(\mu)\circ G(f)^{-1})=\mu \tr\big((f^{1/2}+f^{-1}G(f))\circ
G(f^{-1})\big)\label{ldf0}
\end{eqnarray}
We apply the maximum principle: at a point where  $\mu$ is maximum, $A(\mu)$ is nonpositive, hence $\tr(A(\mu)\circ G(f)^{-1})$ is nonpositive. We also know that
$f^{1/2}+f^{-1}G(f)$ is positive. It follows that  $\tr\big((f^{1/2}+f^{-1}G(f))\circ
G(f^{-1})\big)$ is positive. Hence Equation (\ref{ldf0}) implies that $\mu$ is nonpositive at its maximum. Symmetrically we prove that the minimum of $\mu$ is nonnegative. Hence $\mu=0$.

\qed

\begin{proposition}The operator $D$ is proper.\end{proposition}  
\proof According to the terminology used in the Appendix \ref{ma}, $D$  is a Monge Ampère operator. To prove that $D$ is proper, by Proposition \ref{ma1}, it suffices to find {\em a-priori} bounds  -- depending on $g$ -- for $f$ and its first derivatives whenever $D(f)=g$.

{\em We first obtain $C^0$ bounds on $f$}. Let $k_1=\inf (f)$. At a point where $f$ reaches its minimum, $A(f)$ is a positive operator. It follows that at this point we have
$$
D(f)\leq D(k_1)=k_1^3 D(1).
$$
Therefore
$$
\inf (f)\geq (\frac{\inf(g)}{\sup(D(1))})^{\frac{1}{3}}.
$$
A symmetric argument yields
$$
\sup(f)\leq (\frac{\sup(g)}{\inf(D(1))})^{\frac{1}{3}}.
$$

{\em We now obtain $C^1$ bounds on $f$}:
We restrict  a function $f$ of $\cU$ on any geodesic for $\nabla$. Hence we obtain a function depending on one variable so that
$$
\ddot f \leq f.
$$
Since $f$ is bounded, this implies that $\dot f$ is also bounded.
Therefore $f$ has $C^1$-bounds.
\qed

\begin{proposition} The set \ $\cU$ is connected.
\end{proposition}
  
\proof Indeed $\cU$ is the set of functions $f$ such that
$$
fS_u-\nabla^2f,
$$
is a a positive symmetric tensor. It follows that $\cU$ is convex, hence connected.
\qed

\section{Projective structures and cohomology classes}\label{secOp}

Let $\nabla$ be a flat connection on  a rank 3 vector bundle $E$ over a surface $\Sigma$. We say a 1-form $\alpha$ with values in $E$ is {\em injective} if
$$
\forall X\in T\Sigma\setminus\{0\}, \ \ \alpha(X)\not=0.
$$
Every injective injective closed 1-form $\alpha$ with value in $E$  {\em  defines} a section $v=\alpha(T\sigma)$ of $\mathbb P (E^*)$. Moreover, an injective closed immersive form $\alpha$ defines a symmetric bilinear tensor up to a multiple on $TS$ by
$$
h^\nabla_\alpha(X,Y)=\pi(\nabla_X (\alpha(Y))).
$$
where $\pi$ is the projection of $E$ on the rank-1 bundle $E/\alpha(TS)$. If $h^\nabla_\alpha$ is non degenerate, then $v$ defines an equivariant immersion, hence a projective structure $\mathfrak p$ on $\Sigma$. 

We say that $\omega$ is {\em convex} if  for every non zero vector
$h^\nabla_\omega(X,X)\not=0.$ We say the complex structure $J$ is compatible with the convex form  $\omega$ if it defines the conformal class of $\omega$. From the previous observation a convex 1-form defines a projective structure.

Let $\mathcal O_{\mathfrak p}$  be the  open cone of cohomology classes of  convex  1-form in $H^1(E)$ defining $\mathfrak p$. We shall prove the following theorem
\begin{theorem}\label{Op}
Let  ${\mathfrak p}$ be a projective structure. Let  $J$ be  a complex structure on $S$, then there exists a convex 1-form defining  ${\mathfrak p}$ and compatible with $J$. As a consequence $\mathcal O_{\mathfrak p}$ is non empty.

Finally the projective structure ${\mathfrak p}$ is convex if and only if  $\mathcal O_{\mathfrak p}$ contains $0$ . 
\end{theorem}

\rks \begin{itemize}
\item 
It follows that  if $\mathfrak p$ is convex, then $\mathcal O_{\mathfrak p}=H^1(E)$.
\item This result leads to the following natural set of questions: 
\begin{enumerate}
\item given ${\mathfrak p}$ and $J$ does there exist a better convex form defining ${\mathfrak p}$ and compatible with $J$ ? A positive answer to this question would lead to a parametrisation of $\mathcal O_{\mathfrak p}$ by Teichmüller space. 
\item given ${\mathfrak p}$ and $\omega$ in $\mathcal O_{\mathfrak p}$ does there exist a better complex structure $J$ such that $\omega$ can be represented by a convex form defining ${\mathfrak p}$ and compatible with $J$ ?
\end{enumerate} A positive answer to these questions would lead to a map parametrising $\mathcal O_{\mathfrak p}$ by Teichmüller space.  In Theorem \ref{phodge}, when $\mathfrak p$ is convex,  we actually produce a map from $H^1(E)$ to Teichmüller space. 
\end{itemize}

\subsubsection{A preliminary proposition}
We begin with a  proposition
\begin{proposition}\label{form-rp}
Let $\alpha$ be an immersive convex form on $TS$ with value in $E$ equipped with a flat connection $\nabla$. Then, the geodesics for the associated projective structure are the curves $c$ so that there exist a vector field $Y$ along $c$ so that
$$
\nabla_{\dot c} \alpha(Y) =0.
$$
\end{proposition}
\proof This is a local statement. Therefore, we can assume that the connection $\nabla$ is the trivial one and we identify sections of the vector bundle as maps with values  in a vector space. Then, a curve  $c$ is a geodesic if and only if, its image under $v=\alpha(T\Sigma)$ is a dual projective line. This means that all the planes $\alpha(T\Sigma)$ contain a common non zero vector $Z$. Therefore, if $Y$ is the vector field along $c$ such that $\alpha(Y)=Z$, it follows that $\alpha(Y)$ is a parallel section of $E$ along $c$. \qed

\subsubsection{Proof of Theorem \ref{Op}}
\proof 
Let $J$ be a complex structure on $S$. Let  $\mathfrak p$ be a projective structure on $S$. Let $(\nabla,h)$ be a pair consisting of a  torsion free connection $\nabla$ preserving a volume form $\omega$ and  representing $\mathfrak p$, and $h$ a symmetric tensor so that $\nabla^h$ is flat. Recall that $\nabla^h$ is a connection on $E=TS\oplus \mathbb R$. Let $g$ be the metric  on $TS$ given by $g(X,X)=\omega(X,JX)$. We equip $E$ with the metric given by
$$
G((X,\lambda),(X,\lambda))=g(X,X)+\lambda^2.
$$
Let $\nabla^{h,*}$ (resp. $\nabla^*$) be the dual connection to $\nabla^h$ with respect to $G$ (resp. $g$). Let $H$ be the symmetric tensor so that $g(HX,Y)=h(X,Y)$.
Observe that
$$
\nabla^{h,*}_Z(X,\lambda)=(\nabla^*_Z X +\lambda H (Z), g(X,Z) + d\lambda(Z)).
$$
Let $\alpha\in\Omega^1(TS)\otimes E$ be given by $\alpha(X)=(X,0)$. Observe that $d^{\nabla^{h,*}}\alpha=0$.
Moreover,
$$
h^{\nabla^*}_\alpha=\pi(\nabla^{h,*}X Y)=g(X,Z).
$$
Hence $\alpha$ is convex (with respect to the connection $\nabla^{h,*}$). Let $\mathfrak q$ be the associated projective structure. By Proposition \ref{form-rp}, the geodesics for $\mathfrak q$ are those curves $c$ along which there exist a vector field  $Y$ such that  $\alpha(Y)$ is parallel.  Therefore we have
$$
0=\nabla^{h,*}_{\dot c} (Y,0)=(\nabla^*_{\dot c} Y, g(\dot c,Y))
$$
This is equivalent to the fact that $\nabla^*_{\dot c} J \dot c$ is colinear to $J\dot c$. Since, $\nabla^*=-J\nabla J$, it follows that this is equivalent to the fact that $c$ is a geodesic for $\nabla$. Hence, $\mathfrak q=\mathfrak p$. This finishes the proof of the first part of the Theorem.

Finally, if $0\in \mathcal O_{\mathfrak p}$, this means that we can find a section $u$ so that $\nabla u$ is an immersive form defining  $\mathfrak p$, and $S_u$ is definite positive. Hence $\nabla^u$ defines a convex structure (see the  definitions and notations in Paragraph \ref{connsym})

By Proposition \ref{form-rp}, a curve $c$ is a geodesic for  $\mathfrak p$ if there exists a non zero vector field $Y$ along $c$ so that
$$
\nabla_{\dot c}\nabla_Y u=0.
$$
This means that $\nabla^u_{\dot c} Y=0$ and $S_u(\dot c, Y)=0$. This means that $c$ is a geodesic for the connection $\nabla^{u,*}$ dual to $\nabla^u$ with respect to  $S_u$. 

Finally, we notice that $\nabla^{u,*}$ also defines a convex structure. Indeed, by the previous discussion, the connection $\nabla^*=(\nabla^{u,*})^{S_u}$ is dual to the connection $\nabla=(\nabla^u)^{S_u}$ with respect to the metric
$$
G((X,\lambda),(X,\lambda))=S_u(X,X)+\lambda^2.
$$
Thus $\nabla^*$ is flat and $\nabla^{u,*}$ is convex. We have finished to prove that $\mathfrak p$ is convex. \qed

\section{Cohomology classes and complex structures}

Our aim is to give a description of $H_\rho^1(E)$ is terms of complex structures on the surface when $\rho$ is the holonomy of a convex projective surface.

\begin{theorem}\label{pcohom}{\sc[Complex structures]}
Let $\mathfrak p$ be a convex projective structure of holonomy $\rho$ on a closed surface $\Sigma$.  Let $\nabla$ be the  volume preserving connection and $J_0$ be the complex structure described by Theorem \ref{rp->cub}. Let $\nabla^*=-J_0\nabla J_0$ be the dual connection. Let 
$$
\mathcal J=\{J\in \Gamma({\rm End}(TS))/ J^2=-1 \hbox{ and } d^{\nabla^*}J=0\}.
$$
Then  the map from $\mathcal J$ to $H^1_\rho (E)$ given by
$$
J\to J_0J\in \Gamma({\rm End(TS)})=\Omega^1(S)\otimes TS\subset \Omega^1(S)\times E,$$
is a bijection.
\end{theorem}

Notice that  when $\rho$ is with values in $SO(2,1)$, the canonical map from $
\mathcal J$ to Teichmüller space is an isomorphism. This follows from the fact that  $d^\nabla J=0$ if and only if the identity map is an harmonic mapping from $\Sigma$ equipped with $J$ to $\Sigma$ equipped with $\nabla$.

The Theorem is given by Proposition \ref{cohomj}.

We also prove a linear version of this Theorem which gives another  interpretation of $H^1_\rho(E)$ 
\begin{theorem}\label{phodge}{\sc [Hodge representatives]}
Let $\mathfrak p$ be a convex projective structure of holonomy $\rho$ on a closed surface $\Sigma$.  Let $\nabla$ be the  volume preserving connection and $J_0$ be the complex structure described by Theorem \ref{rp->cub}. Let 
$$
\mathcal H_{J_0}=\{A\in \Gamma({\rm End}(TS))/ AJ_0=-J_0A \hbox{ and } d^{\nabla}A=0\}.
$$
Then  the map from $\mathcal H_{J_{0}}$ to $H^1_\rho (E)$ given by the inclusion
$$
\mathcal H_{J_0}\subset  \Gamma({\rm End(TS)})=\Omega^1(S)\otimes TS\subset \Omega^1(S)\otimes E,$$
is a bijection.
\end{theorem}

The theorem is given by Proposition \ref{phodgeprop}.

We observe that in the case  $\rho$ is with values in $SO(2,1)$, there is a bijection between $\mathcal H_{J_0}$ and the space of quadratic holomorphic differentials and we recover here case the Eichler-Shimura isomorphism.

In the next section, we explain how these results yield  interesting symmetries  between representations in the affine group.

\subsection{Cohomology classes and complex structures}\label{cohomj}
Our aim is to prove the following Proposition which generalises Theorem \ref{pcohom}.

\begin{proposition}
Let $\mathfrak p$ be a convex projective structure on $S$ of holonomy $\rho$. Let $\nabla$ be a connection representing $\mathfrak p$. Assume there exist a metric $g$ such that $\nabla^g$ on $E=TS\oplus\mathbb R$ is flat. Let $J_0$ be the complex structure of $g$. Let $\tilde \nabla =-J_0\nabla J_0$. Let $\mu$ be a an element of $H^1_\rho (E)$. Then there exists a unique complex structure $J$ on $S$ such that
$$
d^{\tilde \nabla} J=0,
$$
and 
$$
J_0J \in \Omega^1 (S)\otimes TS\subset \Omega^1(S)\otimes E,
$$
is in the cohomology class of $\mu$.
\end{proposition}
We observe that in the case $\nabla$ preserves the volume form of $g$ ({\em i.e.} we are in the case described by Theorem \ref{rp->cub}), then $\tilde\nabla$ is the dual connection of $\nabla$ with respect to $g$.

\subsubsection{A Monge-Ampère equation}
Let $g$ be a metric on $S$. Let  $B_\mu$ be any symmetric operator on $TS$. Let $\nabla$ be any connection on $S$.  Let
$$
H_\mu \mapping{ C^\infty(S)}{C^\infty(S)}{f}{\det (B_\mu+ f  - \nabla^2 f)\}}
$$
 
We now prove
\begin{proposition}\label{detbmu} 
Let 
$
\mathcal U_\mu=\{f\in C^\infty(S)/ H_\mu(f) >0\}
$. Then $H_\mu$ is a diffeomorphism form $\mathcal U_\mu$ to $C^\infty (S,]0,\infty[)$.
\end{proposition}
\subsubsection{Proof of Proposition \ref{detbmu}}
This proof  will follow closely the strategy of the proof of Proposition \ref{nu}. We use  the following three steps
\begin{enumerate}
\item $H_\mu$ is a local diffeomorphism from $\mathcal U_\mu$ in $\CI(S,]0,+\infty[)$,
\item $H_\mu$ is proper,
\item $\mathcal U_\mu$ is connected.
\end{enumerate}
This will require three separate propositions.

\begin{proposition}\label{prop:ma3}
$H_\mu$ is a local diffeomorphism from $\mathcal U_\mu$ in $\CI(S,]0,+\infty[)$,
\end{proposition}
  
\proof We first compute the linearised operator $L^{H_\mu}_{f}$ at $f$ of $H_\mu$. Let 
$$G(f)= B_\mu +f- A(f).$$
Then a straightforward computation give
\begin{eqnarray*}
L^{H_\mu}_{f}(g) &=&\frac{1}{H_\mu(f)}\tr\big((g -A(g))\circ G(f)^{-1}\big)
\end{eqnarray*}

Since $f\in\cU_\mu$, $G(f)$ is a positive symmetric operator, it follows that $L^{H_\mu}_{f}$ is an elliptic operator. By the implicit function theorem, to prove the Proposition it suffices to show that $L^{H_\mu}_{f}$ is invertible. Since
$L^{H_\mu}_{f}$ is homotopic to a Laplacian, its index is zero. It thus suffices to show that 
$L^{H_\mu}_{f}$ is injective.
Let $g$ so that  $L^{H_\mu}_{f}(g)=0$. Hence 
\begin{eqnarray}
\tr(A(g)\circ G(f)^{-1})=g\tr\big(
G(f^{-1})\big)\label{ldmu0}
\end{eqnarray}
We apply the maximum principle: at a point where  $g$ is maximum, $A(g)$ is nonpositive, hence $\tr(A(\mu)\circ G(f)^{-1})$ is nonpositive. We also know that
$G(f)$ is positive. It follows that  $\tr\big(
G(f^{-1})\big)$ is positive. Hence Equation (\ref{ldmu0}) implies that $g$ is nonpositive at its maximum. Symmetrically we prove that the minimum of $g$ is nonnegative. Hence $g=0$.

\qed

\begin{proposition}The operator $H_\mu$ is proper.\end{proposition}  
\proof According to the terminology used in the Appendix \ref{ma}, $H_\mu$  is a Monge-Ampère operator. To prove that $H_\mu$ is proper, by Proposition \ref{ma1}, it suffices to find {\em a-priori} bounds  -- depending on $g$ -- for $f$ and its first derivatives whenever $H_\mu(f)=g$.

{\em We first obtain $C^0$ bounds on $f$}. Let $k_1=\inf (f)$. At a point where $f$ reaches its minimum, $A(f)$ is a positive operator. Let in general $\Lambda_\mu$ (resp.
$\lambda_\mu$)  be the greatest (resp. smallest) eigenvalue of $B_\mu$. It follows that at this point we have
$$
H_\mu(f)\leq H_\mu(k_1)\leq (k_1+\Lambda_\mu)^2.
$$
Therefore
$$
\inf (f)\geq \inf(g)^{\frac{1}{2}}-\Lambda_\mu.
$$
A symmetric argument yields
$$
\sup(f)\leq \sup(g)^{\frac{1}{2}}-\lambda_\mu.
$$

{\em We now obtain $C^1$ bounds on $f$}.
if $f\in\cU_\mu$, it follows that the function $f$ restricted to any geodesic satisfies
$$
\ddot f \leq f +h
$$
for some function $f$. Now, $C^0$-bounds on $f$ implies bounds on $\dot f$. It follows that $f$ such that $H_\mu(f)=g$ admits $C^1$-bounds.
 \qed

\begin{proposition} The set \ $\mathcal U_\mu$ is connected.
\end{proposition}
  
\proof Indeed $\mathcal U_\mu$ is the set of functions $f$ such that
$$
B_\mu+ f-\nabla^2f,
$$
is a a positive symmetric tensor. It follows that $\mathcal U_\mu$ is convex, hence connected.
\qed

\subsubsection{Proof of Proposition \ref{cohomj}}

Let $\mu$ be a  a closed 1-form in $\Omega^1(S)\otimes E$. We write
$$
\mu(X)=(B_\mu (X),\alpha_\mu (X))\in TS \oplus \mathbb R=E,
$$
where $B_\mu\in {\rm End}(TS)$ and $\alpha\in T^*(S)$. Let $\xi_\mu$ be the vector field such that
$g(\xi_\mu,X)=\alpha_\mu(X)$. Then 
$$
\alpha_{\mu -\nabla \xi_\mu}=0.
$$
It follows that every cohomology class in $E$ has a representative  $\omega$ such that $\alpha_\omega=0.$

Let $\omega$ be such a representative. Let $v=(\xi,f)$ be a section of $E$. We observe also that 
$$
\alpha_{\omega +\nabla v}=0,
$$
if and only if $\xi=-{\rm grad} f$.
Let $\nabla_*$ be the dual connection to $\nabla$ with respect to $g$, observe that
$$
\nabla {\rm grad} f=\nabla^2_*f.
$$
Let $v_f=(-{\rm grad}f, f)$.
We have
$$
B_{\omega + \nabla v_f}=B_\omega + f- \nabla_*^2f.
$$
By Proposition \ref{detbmu}, we conclude  there exists a unique  representative $\omega$ of the cohomology class $\mu$ such that $\alpha_\omega=0$ and $\det (B_\omega)=1$. 

Proposition \ref{cohomj} follows from Proposition \ref{prop:ma3} and the following two observations
\begin{itemize}
\item if $B\in\Omega^1(S)\otimes TS\subset\Omega^1(S)\otimes E$ is closed, then $d^\nabla B=0$ and $B$ is symmetric.
\item $B$ is symmetric with determinant equal to 1, if and only if $J_0B$ is a complex structure.
\end{itemize}

\subsection{Hodge representatives}

We now prove the first part Theorem \ref{phodge}
\begin{proposition}\label{phodgeprop}
Let $\nabla$ be the  torsion free  connection and $g$ a metric such that $\nabla^g$ is flat of holonomy $\rho$. Let $J_0$ be the complex structure of $g$.  Let 
$$
\mathcal H_{J_0}=\{A\in \Gamma({\rm End}(TS))/ AJ_0=-J_0A \hbox{ and } d^{\nabla}A=0\}.
$$
Then  the map from $\mathcal H_{J_0}$ to $H^1_\rho (E)$ given by the inclusion
$$
\mathcal H_{J_0}\subset  \Gamma({\rm End(TS)})=\Omega^1(S)\otimes TS\subset \Omega^1(S)\otimes E,
$$
is a bijection.
\end{proposition}

\proof  We shall use the notations of the proof of Proposition \ref{cohomj}. Let $\mu$ be a cohomology class in $H^1_\rho(E)$. Using the same approach, we can represent $\mu$ by a 1-form $\omega$ such that $\alpha_\omega=0$. We now remark there exists a unique function $f$ such that
$$
D(f)=\tr (f  -\nabla^2_*f)=-\tr (B_\omega).
$$
Indeed, $D$ is a linear elliptic operator of index $0$ whose kernel is trivial as it is shown by an easy application of the maximum principle.

We also remark that $A$ is symmetric of trace zero if and only if $AJ_0=-J_0A$. Combining these two remarks, we obtain there exist a unique section $v$ of $E$ such that if $\beta=\omega+\nabla v$ then
\begin{eqnarray*}
\alpha_{\beta}&=&0\\
B_\beta J_0&=&-J_0 B_\beta.
\end{eqnarray*}
The statement follows \qed
\section{Dualities and symmetries of moduli spaces}
 We now explain that  Theorems \ref{rp->cub}, \ref{pcohom} and \ref{phodge}, give rise to interesting symmetries of moduli spaces of representations, well known in the first case, but more mysterious in the other cases.
 
\subsection{Contragredient representation and Theorem \ref{rp->cub}}
 
We define ${\rm Rep_H}(\grf,SL(3,\mathbb R))$ to be the component of  the space of representations which contains the cocompact representations in $SO(2,1)$. Let  $SA\!f\!\!f (3,\mathbb R)$ be the special affine group in dimension 3 
$$
SA\!f\!\!f (3,\mathbb R)=\mathbb R^3\rtimes SL(3,\mathbb R).
$$
We define similarly ${\rm Rep_H}(\grf, SA\!f\!\!f (3,\mathbb R))$ to be the set of those representations whose linear part is in ${\rm Rep_H}(\grf,SL(3,\mathbb R))$. We observe that by results of Choi and Goldman, ${\rm Rep_H}(\grf,SL(3,\mathbb R))$  is precisely the set of monodromies of convex projective structures.

Let $\omega$ be a volume form on $S$. By our Theorem \ref{rp->cub}, it follows that 
$${\rm Rep_H}(\grf,SL(3,\mathbb R)
$$ 
is in one-to-one correspondence with the space of triples $(\nabla,\omega,J)$ where $\nabla$ is  connection, $\omega$ is a volume form, and $J$ is a complex structure, which satisfy condition (H):
\begin{eqnarray}
\left\{
\begin{array}{rcl}
\nabla_XY-\nabla_YX&=&[X,Y]\\
\nabla\omega&=&0\\
d^\nabla J&=&0\\
R^\nabla(X,Y)JZ&=&\omega(X,Y)Z.\\
\end{array}
\right.\label{condH}
\end{eqnarray}
It is a trivial observation that $(\nabla,\omega,J)$ satisfies condition (H), if and only if $(-J\nabla J,\omega,J)$ does. The corresponding duality in the space of representations is the the duality which associates to a representation its contragredient representation as it is shown by an easy exercise left to the reader. Its set of fixed point is the space of representations with monodromy in $SO(2,1)$.

\subsection{Theorem \ref{pcohom}  and an involution on the moduli space of representations in the affine group}  

Theorem \ref{pcohom}  provides a more mysterious duality. By this result, 
$${\rm Rep_H}(\grf, SA\!f\!\!f(3,\mathbb R))$$  is in bijection with space of quadruples $(\nabla, \omega,J,J_1)$ where $\nabla$ is a connection $\omega$ a volume form, and $J$ as well as $J_1$ are complex structures, which satisfy condition (I):
\begin{eqnarray}
\left\{
\begin{array}{rcl}
\nabla_XY-\nabla_YX&=&[X,Y]\\
\nabla\omega&=&0\\
d^\nabla J&=&0\\
R^\nabla(X,Y)JZ&=&\omega(X,Y)Z\\
d^\nabla JJ_1&=&0
\end{array}
\right.\label{condI}
\end{eqnarray}
An exercise shows that $(\nabla,\omega,J,J_1)$ satisfies condition (I) if and only if 
$$(-J_1\nabla J_1,\omega,-J_1JJ_1,J_1)$$
 does. We obtain therefore a duality on ${\rm Rep_H}(\grf, SA\!f\!\!f(3,\mathbb R))$ which extends the duality on ${\rm Rep_H}(\grf, SL(3,\mathbb R))$ considered as a subset. However this duality does not respect the projection from ${\rm Rep_H}(\grf, SA\!f\!\!f(3,\mathbb R))$ to ${\rm Rep_H}(\grf, SL(3,\mathbb R))$ and does not seem to have an algebraic description. Again, its set of fixed point is the space of representations with monodromy in $SO(2,1)$.

\subsection{Theorem \ref{pcohom}  and a fourth order symmetry  on the moduli space of representations in the affine group}  
Finally, Theorem \ref{phodge} also provides a symmetry, of order 4, on 
$${\rm Rep_H}(\grf, SA\!f\!\!f(3,\mathbb R)).$$ 
By this result, this moduli space  is in bijection with space of quadruples  
$(\nabla,\omega, J,A)$ where $\nabla$ is a connection, $\omega$ a volume form, $J$ is complex structure and $A$ is an endomorphism of $TS$, which satisfy condition (J):
\begin{eqnarray*}
\left\{
\begin{array}{rcl}
\nabla_XY-\nabla_YX&=&[X,Y]\\
\nabla\omega&=&0\\
d^\nabla J&=&0\\
R^\nabla(X,Y)JZ&=&\omega(X,Y)Z.\\
d^\nabla A&=&0\\
AJ&=&-JA
\end{array}
\right.
\end{eqnarray*}
Again, it is an exercise that $(\nabla,J,A)$  satisfy condition (J) if and only if $(-J\nabla J,J, JA)$ does.  Observe that the map
$$
j:(\nabla,J,A)\to (-J\nabla J,J, JA),
$$
is actually of order 4: $j^2$ sends  $(\nabla,J,A)$ to $(\nabla,J,-A)$ and is the antipody on the vector bundle
$${\rm Rep_H}(\grf, SA\!f\!\!f(3,\mathbb R))\rightarrow{\rm Rep_H}(\grf, SL(3,\mathbb R)).$$ 
This mapping $j$ extends the duality of $Rep_H(\grf,SL(3,\mathbb R)$ as a subset and also factors over the projection to  $Rep_H(\grf,SL(3,\mathbb R)$.

\section{An affine differential  interpretation}\label{sec:aff}

In this section we give an interpretation of our Theorems \ref{rp->cub}, \ref{pcohom} in terms of affine differential geometry. We also give an interpretation of Theorem \ref{Op} in this language.

We begin by recalling briefly the description of convex hypersurfaces in affine differential geometry. Let $E$ be an affine space equipped with a constant volume form $\Omega$. We denote by $D$ its connection. Let $\Sigma$ be a locally convex hypersurface in $E$. Then there exists a unique pair $(g,\nu)$ such that 
\begin{itemize}
\item $\nu$ is a vector field along $\Sigma$ transverse to $T\Sigma$, 
\item $g$ is a metric on $\Sigma$ whose volume form is $i_\nu \omega$,
\item for all $X$ in $ T\Sigma$, we have  $\nabla_X\nu \in T\Sigma$,
\item for all $X,Y$ in vector fields on $\Sigma$ we have $D_X Y -g(X,Y)\nu \in T\Sigma$.
\end{itemize}
The vector field is the {\em affine normal vector field} and $g$ is the {\em Blaschke metric}. We call $B=\nabla\nu$ the {\em affine shape operator} and $\nabla$ is the Blaschke connection

In other words, if  we decompose 
$$TE\big\vert_\Sigma=T\Sigma\oplus \mathbb R\nu,$$
Then the the connection $D$ on $\Sigma$ decomposes as
\begin{eqnarray}
D_X(Y,\lambda)=(\nabla_XY + \lambda B(X), g(X,Y) +d\lambda(X)).\label{affdef}
\end{eqnarray}

If $\bar\nabla$ is the Levi-Civita connection of the metric $g$, $\nabla-\nabla_g$ is the {\em Pick invariant} $P$. Conversely, if $\Sigma$ is simply connected, if $B$, $\nabla$ and  $g$  satisfy the following conditions
\begin{itemize}
\item $\nabla$  is torsion free, preserves the volume form of $g$ 
\item the connection $D$ defined by Equation (\ref{affdef}) is flat, that is
\begin{eqnarray}
\left\{
\begin{array}{rcl}
  g(BX,Y)&=&g(X,BY)\\
 d^\nabla g&=&0\\
 d^\nabla B&=&0\\
R^\nabla(X,Y)Z&=&g(X,Z)BY-g(Y,Z)BX
\end{array}
\right.\label{condflat}
\end{eqnarray}
\end{itemize}
then there exists an immersion of  $\Sigma$ whose Blaschke metric is $g$, shape operator is $B$ and Blaschke connection is $P$. If $\Sigma$ is not simply connected, then the universal cover of $\Sigma$ possesses 
an immersion equivariant under a representation in the special affine group. The linear part of this representation is given by the holonomy ${\rho}_{0}$ of $D$, and the affine extension as en element of $H^1_{\rho_{0}}(\pi_{1}(\Sigma),E)=H^1_{D}(T\Sigma\oplus\mathbb R)$ is represented  by the element $\omega$ of $\Omega^1(\Sigma)\otimes (T\Sigma\oplus\mathbb R)$ given by $\omega(X)=(X,0)$.

We concentrate now on the case of surfaces and write every metric $g$ as $g=\omega(.,J)$, where $\omega$ is the volume form of $g$ and $J$ is the complex structure of $g$. 
We write now the equations on $\nabla$, $\omega$, $J$ and $B$ which translates the condition (\ref{condflat}) above
\begin{eqnarray}
\left\{
\begin{array}{rcl}
 \nabla_{X}Y-\nabla_{Y}X&=&[X,Y]\cr
 \nabla\omega&=&0\cr
 \tr (BJ)&=&0\cr
 d^\nabla J&=&0\cr
 d^\nabla B&=&0\cr
R^\nabla(X,Y)Z&=&-\omega(X,Y)BJZ
\end{array}
\right.\label{condB}
\end{eqnarray}

\subsection{Hyperbolic affine spheres and Theorem  \ref{rp->cub} }

An {\em affine sphere} is such that $B=k{\rm Id}$. For $k=-1$, we say the affine sphere is {\em elliptic}, for $k=0$ we say it is {\em parabolic}, for $k=1$, we say  it is {\em hyperbolic}. There is a strong relations between strictly convex cones and hyperbolic affine spheres as is explained by the following difficult theorem conjectured by Calabi \cite{Ca}. This result is due to Cheng and Yau \cite{CY1} and  \cite{CY2} later completed and clarified by the work of Gigena \cite{Gi}, Sasaki \cite{Sa} and A.M.\ Li \cite{Li}, \cite{Li2}.
\begin{theorem}{\sc [Cheng-Yau] [hyperbolic affine spheres bound cones]}
If $\Sigma$ is a hyperbolic affine sphere with a complete Blaschke metric, then $\Sigma$ is properly embedded and bounds a  convex cone. Conversely any strictly convex cone is asymptotic to such a unique hyperbolic affine sphere.
\end{theorem}

In the case where $E$ is of dimension 3, it follows from this result that for every convex structure one can associate an hyperbolic affine sphere invariant under the monodromy of the convex structure: the affine sphere asymptotic to the convex set of $\mathbb P (E)$ on which the monodromy acts cocompactly. This is precisely the content of Theorem \ref{rp->cub}. However, our proof uses a simpler approach.

The relation of Theorem \ref{cub->rp}, between cubic holomorphic form and the Pick invariant is a Theorem of C.-P. Wang \cite{wcp} as we already said.

\subsection{Constant Gaussian curvature surfaces and  Theorem \ref{pcohom}}
We suppose again that $E$ is of dimension 3.
We consider now {\em constant Gaussian curvature 1 affine hypersurfaces} (or CCG hypersurfaces) namely those surfaces for which $\det B=1$ where $B$ is the affine shape operator. Therefore, we can write $B=JJ_{1}$ where $J_{1}$ is a complex structure on $\Sigma$.

It follows from an easy check that $(\nabla,\omega,J,B=JJ_{1})$ satisfies Condition \ref{condB}, if and only if 
$(-J\nabla J, J_{1},J)$ satisfies Condition \ref{condI}. In other words, we can restate Theorem \ref{pcohom} using this observation in the following way.
\begin{theorem}
Let $S$ be a compact surface. Given any representation $\rho$ element of 
${\rm rep_{H}}(\grf, SA\!f\!\!f(3,\mathbb R)$, there exists a unique CCG invariant under this representation.
\end{theorem}
This result does not seem to be known in the affine differential world.

\subsection{Interpretation of Theorem \ref{Op}}  Every locally convex surface $S$ in the affine space admits a natural projective structure: the one given by the immersion $s\mapsto T_{s}S$. Whenever $S$ is equivariant under a representation $\rho$, the holonomy $\dot\rho$ of  flat connection $D$ on $TS\oplus \mathbb R$ described by Formula \ref{affdef} is the linear part of $\rho$. Moreover, the element $\omega$of $\Omega^1(TS)\otimes (TS\oplus\mathbb R)$ defined by $\omega(X)=(X,0)$ is a representative in $H^1_{\rho_{0}}(\mathbb R^3)$ of the cohomology class describing the extension from $\rho_{0}$ to $\rho$. Moreover $\omega$ is convex in the sense of Section \ref{secOp}. Conversely, every convex closed 1-form is obtained this way.

Therefore we can reinterpret Theorem \ref{Op} in the following way.

\begin{theorem}
Given a projective structure $J$ and a complex structure $\mathfrak p$ on $S$, there exists a locally convex surface $S$ in the affine three dimension space whose Blaschke metric is conformal to $J$ and which defines $\mathfrak p$ equivariant under a representation $\rho$ whose linear part is the monodromy $\rho_{0}$of the projective structure $\mathfrak p$. We can choose $\rho$ to be be conjugate to $\rho_{0}$ if and only if $\mathfrak p$ is convex.
\end{theorem}

\section{A Higgs bundle interpretation}

We now recall briefly the work on Hitchin on representations of surface groups in $PSL(n,\mathbb R)$ and explain using \cite{FLEnergy} how it fits with the present work. 

\subsection{Representations and harmonic mappings}

Following \cite{FLA}, we define a {\em Fuchsian representation} of $\grf$ in $PSL(n,\mathbb R)$ to be a representation which factors through the irreducible representation of $PSL(2,\mathbb R)$ in  $PSL(n,\mathbb R)$ and a cocompact representation of $\grf$ in $PSL(2,\mathbb R)$. A {\em Hitchin representation} is a representation which may be deformed in a Fuchsian representation. The space of Hitchin representation is denoted by 
$${\rm Rep}_H(\grf,SL(n,\mathbb R))
$$ 
and is called a Hitchin component.

In his article  \cite{H}, N. Hitchin gives  explicit parametrisations of Hitchin
components. Namely, given a choice of a complex structure
$J$ over  a given compact surface $S$, he produces a homeomorphism
$$
H_J: \mathcal Q(2,J)\oplus\ldots\oplus\mathcal Q(n,J) \rightarrow{\rm Rep}_H(\grf,SL(n,\mathbb R)),
$$
where $\mathcal Q(p,J)$ denotes the space of holomorphic
$p$-differentials on the Riemann surface $(S,J)$.  The main
idea in the proof is  first to identify representations with
harmonic mappings as in K. Corlette's seminal paper
\cite{KC}, (see also  \cite{FL1}), second to use the fact  a harmonic
mapping $f$ taking values in a symmetric space gives rise to
holomorphic differentials in manner
similar to that in which  a connection gives rise to differential forms
in Chern-Weil theory (cf. Paragraph 7.1.2 of \cite{FLEnergy}). We explain quickly the construction. Namely, we have a parallel symmetric $p$ form $q_p$ on $M=SL(n,\mathbb R)/SO(n,\mathbb R)$. Identifying (after a choice of a base point) $M$ with the space of metrics of volume 1 and $T_g M$ with the space of self-adjoint (with respect to $g$),  endomorphism of $\mathbb R^n$ we set
$$
(q_p)_g (A,\ldots,A)=\tr (\underbrace{A \ldots A}_{p \hbox{  \tiny times}}).
$$
Then we can complexify $q_p$ as a parallel symmetric complex $p$ form on the complexified tangent bundle. Then, the $p$-ic holomorphic form $Q_p(f)$ associated to a harmonic mapping $f$ with values in $M=SL(n,\mathbb R)/SO(n,\mathbb R)$ is
$$
Q_p(f)=q_p (\underbrace{T_{\mathbb C}f \ldots T_{\mathbb C} f}_{p \hbox{  \tiny times}}),
$$
where $T_{\mathbb C}f$ is the complexification of $Tf$ : 
$$
T_{\mathbb C}f (u)=Tf (u)-iTf (Ju).
$$
We observe that $Q_{2}(f)=0$ if and only if $f$ is minimal ({\it cf} Proposition 7.1.3 of \cite{FLEnergy} and \cite{SU2},\cite{SY}).

We also observe that  given a Hitchin representation $\rho$ and a complex structure  $J$ on $S$, on obtain a number $e_{\rho}(J)$, the energy of the associated harmonic mapping. The {\em energy} will be the function on Teichmüller space given by $J\mapsto e_{\rho}(J)$. 

\subsection{Representations, energy  and minimal surfaces}
However one drawback of this construction is that $H_J$ depends on the choice of the complex structure $J$. In particular, it breaks the invariance by the Mapping Class Group and therefore this construction does not give information on the topological nature of  ${\rm Rep}_H(\grf,SL(n,\mathbb R))/\mathcal M
  (S)$. We explain now a more equivariant (with respect to the action of the Mapping Class Group) construction. Let 
$\mathcal E^{(n)}$ be the vector bundle over Teichmüller space whose fibre above the complex structure $J$ is
$$
\mathcal E^{(n)}_J=\mathcal Q(3,J)\oplus\ldots\oplus\mathcal Q(n,J).
$$
We observe that the dimension of the total space of $\mathcal E^{(n)}$ is the same as that of ${\rm Rep}_H(\grf,SL(n,\mathbb R))$ since the dimension of the "missing" quadratic differentials in $
\mathcal E^{(n)}_J$ accounts for the dimension of Teichmüller space.
account for .
We now define the {\em Hitchin map} 
$$
H \mapping{\mathcal E^{(n)}}{{\rm Rep}_H(\grf,SL(n,\mathbb R)}{(J,\omega)}{H_J(\omega).}
$$
We are aware that this terminology is awkward since this Hitchin map  is some kind of an inverse of what is usually called the Hitchin fibration.
From Hitchin construction, it now  follows this map is equivariant with respect to the Mapping Class Group action. We quote from \cite{FLEnergy} the following two results
\begin{theorem}\label{hitchproper}
The energy $e_{\rho}$ on Teichmüller space is proper.
\end{theorem}

\begin{theorem}\label{hitchsurjint}
The Hitchin map is surjective.
\end{theorem}

Our strategy is to identify $\mathcal E^{(n)}$ with the moduli space of  equivariant minimal surfaces  in the associated symmetric space and to prove that there exists an equivariant minimal surface for every representation by tracking a critical point of the energy. Indeed, harmonic mappings for which the quadratic differential vanishes are conformal, and minimal surfaces are critical points of the energy (\cite{FLEnergy} and \cite{SU2},\cite{SY}).

Our conjecture in \cite{FLA} is that the Hitchin map is a homeomorphism, which is also equivalent by the above discussion to the following one  
\begin{conj}
If $\rho$ is a Hitchin representation, then $e_\rho$ has a unique  critical point. 
\end{conj}

\subsection{The case of $n=3$.}

Let $S$ be a locally convex surface in a three affine space $E$ equipped with a volume form. Then we define the {\em Blaschke lift} $G$ as a map from $S$ to the space ${\rm Met}(E)$ of euclidean metrics on $E$ of volume 1:
$$
\mapping{S}{{\rm Met}(E)}{s}{G(s)\hbox{ such that } 
G(s) (X,\lambda)=g_s(X,X) + \lambda^2,}
$$
where $E$ is identified with $T_sS\oplus R$ and $g_s$ is the Blaschke metric on $T_sS$.

We know prove the following proposition whose first part is well known.
\begin{proposition}
$S$ is an affine sphere if and only if $G(S)$ is a minimal surface. Moreover, real part of the cubic holomorphic $q_3(G)$  associated to $G$ is related to the Pick invariant $A$ of $S$ : 
\begin{eqnarray}
12\cdot g_s(A(X)Y,Z)&=&\Re ( Q_3(G)(X,Y,Z)).
\end{eqnarray}
\end{proposition}
\proof We just prove the second part of the proposition. We first observe that by definition identifying $T_gM$ with the space of symmetric endomorphisms of $\mathbb R^3$. We denote by $X^*$ the transpose of $X$,
$$
TG (X)=\Gamma(X)=
\left(
\begin{array}{cc}
A(X)&X\\ X^*&0
\end{array}
\right).
$$
Then
$$
\Gamma(X)\Gamma(Y)=
\left(
\begin{array}{cc}
A(X)A(Y)+XY^*&A(X)Y\\ X^*A(Y)&X^*Y
\end{array}
\right),
$$
and
\begin{eqnarray*}
\tr(\Gamma(Z)\Gamma(X)\Gamma(Y))&=&
\tr(
A(Z)A(X)A(Y))\\
& &+g(A(Z)X,Y)+g(Z,A(Y)X)+g(Z,A(X),Y)
\end{eqnarray*}
If $S$ is an affine sphere then, $g(A(Z)X,Y)$ is symmetric in $X,Y,Z$ and 
\begin{eqnarray}
A(Z)J=-JA(Z)\label{AJJA}
\end{eqnarray}
It follows that
\begin{eqnarray*}
\tr(\Gamma(Z)\Gamma(X)\Gamma(Y))&=&
3\cdot g(A(X)Y,Z).
\end{eqnarray*}
Now, using the definition of $T_{\mathbb C}G$ we get
\begin{eqnarray*}
\Re(Q_3(G)(X,Y,Z))&=&\Re(\tr (T_{\mathbb C}G(X),T_{\mathbb C}G(Y),T_{\mathbb C}G(Z)))\\
&=&\tr (\Gamma(X)\Gamma(Y)\Gamma(Z)) -\tr (\Gamma(X)\Gamma(JY)\Gamma(JZ))\\
&-&\tr (\Gamma(JX)\Gamma(Y)\Gamma(JZ))-\tr (\Gamma(JX)\Gamma(JY)\Gamma(Z))
\end{eqnarray*}
We also observe that by Equation (\ref{AJJA}) $$
g(A(X)JY,JZ)=-g(A(X)Y,Z).$$
Hence
\begin{eqnarray*}
\Re(Q_3(G)(X,Y,Z))&=&12\cdot g(A(X)Y,Z)\end{eqnarray*}\qed

As an immediate corollary, we obtain

\begin{theorem}
For $n=3$, the Hitchin map is a diffeomorphism. Moreover, the energy $e_{\rho}$ has a unique critical point on Teichmüller space which is an absolute minimum.
\end{theorem}

Indeed, by the previous proposition, the map which associate to a  Hitchin representation the Pick invariant of the associated affine sphere and its complex structure is the inverse (up to normalisation by $1/12$) of the Hitchin map. The second part follows from the fact that a complex structure is a critical point of the energy if and only if the associated harmonic map is minimal.

Therefore, since a Hitchin representation is discrete and torsion free ({\em cf } \cite{FLA}), we also have the following corollary

\begin{coro}\label{cor:min}
Let $\rho$ be a Hitchin representation of $\grf$ in $SL(3,\mathbb R)$. Then, there exists a unique minimal surface $S$ in  $\rho(\grf)\backslash SL(3,\mathbb R)/SO(3,\mathbb R)$ such that the injection is a homotopically equivalence.
\end{coro}
\section{A holomorphic interpretation}\label{sec:holo}
We finish this paper by another interpretation of Condition (E) \ref{condE}. We consider the homogeneous space $M=SL(3,\mathbb R)/SL(2,\mathbb R)$, where $SL(2,\mathbb R)$ is embedded reducibly in $SL(3,\mathbb R)$. The space $M$ is the space of pairs $(P,u)$ such that $P$ is a plane in $\mathbb R^3$ and $u$ is a transverse vector to $P$. 

We observe that  we have the following identification
\begin{eqnarray*}
T_{(P,u)}&=& {\rm Hom}(P,\mathbb Ru)\oplus {\rm Hom}(\mathbb Ru,P)\oplus {\rm Hom}(\mathbb Ru, \mathbb Ru,)\\
&=&\underbrace{P^*\oplus P}_{W}\oplus \ \mathbb R.
\end{eqnarray*}
We identify $P^*$ with $P$ using the 2-fom $i_u \Omega$, where $\Omega$ is the volume form of $\mathbb R^3$. We equip $W$ with the complex structure
$$
J(u,v)=(-v,u).
$$
Our last interpretation is the following.

\begin{theorem}  Given a Hitchin representation $\rho$ of $\grf$  in $SL(3,\mathbb R)$, there exists a unique surface $\Sigma$ everywhere tangent to $W$, complex and equivariant under $\rho$.
\end{theorem}
\proof  Indeed, if $S$ is equipped with a convex projective structure, we obtain map $i$ from $S$ to $M$ from the splitting of Theorem \ref{theo:split}. Then, one checks easily that $i(S)$ is tangent to $W$ and complex. For more details, see Paragraph 2.6 of \cite{maudin}. \qed

\section{Appendix A: geodesics}
We prove the following lemma. Let $M$ be a compact manifold. Let $\nabla$ be a connection on $M$. Let $\gamma$ be a geodesic defined on a maximal interval $I=]a,b[$.  Let $g$ be an auxiliary  metric.  We prove
\begin{lemma}\label{complete}
Let $\mu(t)=\sqrt{g(\dot\gamma(t),\dot\gamma(t))}$. 
Then, for any $c\in I$
\begin{eqnarray}
\int_c^b \mu(s)ds&=&+\infty.\label{infinitelength}
\end{eqnarray}
Moreover, there exists a constant $K$ such that
\begin{eqnarray}
\vert \dt{\mu}\vert \leq K\mu^2.\label{boundedgrwoth}
\end{eqnarray}
\end{lemma}
\proof We now that $\nabla g$ is bounded since $M$ is compact. It follows that
\begin{eqnarray}
\vert \dt{\mu^2}\vert =\vert \nabla_{\dot{\gamma}(t)}g(\dot{\gamma}(t), \dot{\gamma}(t))\vert\leq K\mu^3.
\end{eqnarray}
The second assertion follows. It follows that if $t>s$ we have
\begin{eqnarray}
-K(t-s)\leq \frac{1}{\mu(t)}-\frac{1}{\mu(s)}\leq K(t-s).\label{ineqcomplete}
\end{eqnarray}
We split the end of the proof in two cases
\begin{enumerate}
\item $b=\infty$. Therefore, we have from (\ref{ineqcomplete})
$$
\forall t>c,\  \mu(t)\geq \frac{\mu(c)}{K\mu(c)(t-c)+1}.
$$
Hence
$$
\int_c^\infty \mu(s)ds\geq \int_c^\infty \frac{\mu(c}{K\mu(c)(s-c)+1}ds=+\infty.
$$
\item $b<\infty$. Then, since the geodesic is maximal, it leaves every compact set  in the tangent bundle :
$$
\lim_{t\rightarrow b}\mu(t)=\infty.
$$
The Inequality (\ref{ineqcomplete}) yields
$$
 -K(b-s)\leq-\frac{1}{\mu(s)}
$$
Hence 
$$
\int_c^b\mu(s)ds \geq \int_c^b\frac{1}{b-s}ds=\infty
$$ and the assertion follows.
\end{enumerate}
\qed

\section{Appendix B: Elliptic Monge-Ampère equations.} \label{ma}

Let $S$ be a closed surface.  Let 
 $J^{1}(S,\RR)$ be the space of 1-jets of  functions on $S$ We denote by $j^1f(x)$ the 1-jet of the function $f$ at $x$. Let
 \begin{itemize}
 \item $G$ be a metric on $S$ and $\nabla$ its whose Levi-Civita connexion,
 \item  $d$ be an $\mathbb R^+$-valued function on 
$J^{1}(S,\RR)$,
 \item $W$ a mapping from
 $J^{1}(S,\RR)$ to ${\rm End}(TS)$  such that $W(j_1f(x))$ is a symmetric endomorphism of  $T_{x}(S)$.
\end{itemize}

By definition the Monge-Ampère operator associated to  $G,d,W$ is
par
$$
M(f)=d(j^{1}f)\det(\nabla^{2}f+W(j^1f))
$$
where $\nabla^{2}f$ is defined by:
$$
L_X\cdot L_Y\cdot f-L_{\nabla_{X}Y}\cdot f=g(\nabla^{2}f(X),Y)\, .
$$
We will show, using holomorphic curves, the following proposition
\begin{proposition}\label{ma1} Let $\{g_{n}\}\inn$ be a sequence of positive functions converging $\CI$ to a positive function $g_{0}$. Let  $\{f_{n}\}\inn$ be a sequence of functions such that
\begin{itemize}
\item there exists  a constant $K$ such that, for all $n$, 
$\Vert f_{n}\Vert_{C^{1}}\leq K$
\item $M(f_n)=g_n$.
\end{itemize}
 Then after extracting a subsequence, $\{f_{n}\}\inn$ converges 
 $\CI$ to 
$f_{0}$ such that
$$
M(f_{0})=g_{0} .
$$
\end{proposition} 
We remark that  a general theory of Monge-Ampère geometries related to holomorphic curves has been developed in 
\cite{labgafa}.

\proof We first show that given a positive function $g$, there exists a complex structure $J_g$ defined on the contact  subbundle $P$ of the tangent bundle of  $J'(S,\RR)$,  such that if 
for all $f$ such that  $M(f)=g$, then the graph  $J^{1}(f)$ is holomorphic with respect to $J_g$.

We begin by describing the geometry of 
$J'(S,\RR)=T^*S\times\RR$. The connexion $\nabla$ gives rise to a decomposition of $TJ^{1}(S,\RR)$
$$
TJ^{1}(S,\RR)=TS\oplus T^{*}S\oplus\RR.
$$
In this decomposition, if we see  $j^{1}f$ as a mapping from 
$S$ to $J'(S,\RR)$, then
$$
Tj^{1}f(u)=(u,\nabla_{u}df,df(u)).
$$
The contact subbundle  $P$ at a point
 $(\omega,\la)$ of 
$J'(S,\RR)$ is 
$$
P_{(\omega,\la)}=\{(u,\a,\omega(u)\, ;\quad u\in TS\, , \quad \a\in T^{*}S\}.
$$
Let $i$ be the isomorphism of $TS$ with its dual coming from the metric $G$. Let $J_0$ the complex structure on $TS$ given by the metric. We identify  $P$ with $TS\oplus TS$ using the following isomorphism 
$$
\psi_g :\mapping
{TS\oplus TS}
{P_{(\omega,\la)}}
{(u,v)}
{\big(u,i(-\sqrt{gd^{-1}}J_{0}v -W(u)), \omega(u)\big)}
$$

Let $J_{1}$ be the complex structure on  $TS\oplus TS$ given by~:
$J_{1}(u,v)=(-v,u)$. Let 
$$
J_{g}=\psi_g\circ J_{1}\circ \psi^{1}_g.$$

We now show that  $M(f)=g$, if and only if 
 $T(j^{1}f)(TS)$ is stable by  $J_{g}$. 
 
 This relies on the following observation: if $B$ is a symmetric positive operator on $TS$, then 
 $$(J_{0}B)^{2}=-(\det B)^{2}.
 $$. In our case, we apply this observation to
$$
B(u)=\sqrt{dg^{-1}} \ \ (\nabla_{u}\nabla f+W(u))\, .
$$

Finally, we equip 
$TS\oplus TS$ - and hence  $P$ - with the product metric  $G_{0}$. 
We now show :
$$
{\rm Area}(j^{1}f(S))\leq A(\|f\|_{C^{1}})\, + C.
$$
Indeed
\begin{eqnarray}
{\rm Area}(j^{1}f(S)) &=&\int_{S}\det(1+B)\\ 
 &=&{\rm Area}(S)+\int_{S}\tr(B) + \int_{S}\det(B).
 \end{eqnarray}
Since $M(f)=g$, we have $\det(B)=1$. Moreover, $B$ being a positive operator
\begin{eqnarray}
\int_{S}\tr(B) &\leq& \sup(\sqrt{dg^{-1}})\big(\int_{S}\D f+\int_{S}\tr
(W)\big)\\ 
 &\leq& C(\|f\|_{C_{1}}) +D,
 \end{eqnarray}
where the constants $C$ and $D$ only depend on $g$ and $M$.

We can now translate the hypothesis of our Proposition in an holomorphic language.
It follows from our construction that the sequence of graphs $\{S_n\}\inn$
$\{j^{1}(f_{n})\}\inn$, where  $f_{n}$ satisfies $M(f_{n})=g_{n}$, is a sequence of holomorphic curves of bounded area
for a sequence of converging complex structures $\{J_{g_{n}}\}\inn$. 
As it is described in \cite{maudin}, we can apply to this situation Gromov's compactness theorem. Thus the sequence  $\{S_n\}\inn$  converges -- after taking a subsequence -- to a holomorphic curve moduli the apparition of bubbles.

For topological reasons, since $S_n$ are graphs over $S$, the bubbles are subset of the fibres. Therefore, no bubbles can occur since the tangent space of fibre does not contain any complex subspace. It follows that our sequence of graphs converges smoothly to a graph. Hence that $\{f_n\}$ converges $\CI$. \qed

\section{Appendix C:  Laplace equations}\label{lapla}

Let as before  $S$ be a closed surface.  Let 
 $J^{1}(S,\RR)$ be the space of 1-jets of  functions on $S$.  Let
 \begin{itemize}
 \item $G$ be a metric on $S$ and $\nabla$ its whose Levi-Civita connexion,
 \item  $F$ be an $\mathbb R$-valued function on 
$J^{1}(S,\RR)$ 
\end{itemize}
The associated Laplace equation is
$$
L_{F}(f)=\D f +F(j^{1}f)=0.
$$
We are going to sketch -- using holomorphic curves --- a proof of the following classical result

\begin{proposition} Let  $\{F_n\}\inn$ be a sequence of functions converging smoothly 
to $F_{0}$. Let $\{f_{n}\}\inn$ be a sequence  $C^{1}$-bounded functions
such that  
\begin{itemize}
\item $
\exists\, K/\ \forall n, \ |f_{n}\|_{C^{1}}\leq K,
$
\item $
L_{F_{n}}(f_{n})=0\, .
$
\end{itemize}
Then after extracting a subsequence, $\{f_{n}\}\inn$ converges  $\CI$ to a function 
$f_{0}$ such that
$
L_{F_{0}}(f_{0})=0.
$
\end{proposition}

\proof The proof being very similar to the previous one, we are going to be sketchy. In \cite{maudin}, we showed there exists a complex structure $J_F$  on the contact subbundle $P$, such that the graph of $j^1f$ is holomorphic if and only if  $L_F(f)=0$. 

WE now have to obtain a control of the area. 
\begin{eqnarray}
{\rm Area}(j^1f(S)) &\leq& \int_{S}\sqrt{\det (1+(\nabla^{2}f)^{2}})\\ 
&\leq& \int_{S}1+\tr((\nabla^{2}f)^{2})
\end{eqnarray}
Moreover
\begin{eqnarray}
\int_{S}\tr((\nabla^{2}f)^{2}) &=&\int_{S}\sum_{i}\langle \nabla_{X_{i}}\nabla f
| \nabla_{X_{i}}\nabla f\rangle\\ 
&=&-\int_{S}\sum\langle \nabla_{X_{i}}\nabla _{X_{i}}\nabla f
| \nabla f\rangle.
\end{eqnarray}
It is quite classical to show that 
$$
\nabla \D f - \nabla _{X_{i}}\nabla_{X_{i}}\nabla f =H(j^1f)
$$
just depends on the 1-jet of $f$. In particular,
\begin{eqnarray}
\int_{S}\tr((\nabla^{2} f)^{2}) &\leq& C_{1}(\|f\|)_{C_{1}}- \int_{S}\langle
\nabla\D f |\nabla f\rangle\\  
&\leq& C_{2}(\|f\|_{C_{1}})+\int(\D f)^{2}\\ 
&\leq& C_{3}\|f\|_{C_{1}}+C_4.
\end{eqnarray}
Therefore, as before, $j^1f_n(S)$ is a sequence of holomorphic curves of bounded area for a sequence of converging complex structures. Thus , we have convergence up to apparition of bubbles in the fibre. However, in this case the fibre does not contain any compact holomorphic curves. We therefore that $\{f_n\}\inn$ converges after extracting a subsequence.

\auteur
\end{document}